\documentclass[12pt]{article}
\usepackage[utf8]{inputenc}

\usepackage{graphicx}
\usepackage{fullpage}
\usepackage{amssymb}
\usepackage{amsmath}
\usepackage{float}
\usepackage{enumerate}
\usepackage{booktabs}
\usepackage{geometry}
\usepackage{subcaption}
\usepackage{listings}
\usepackage{xcolor}
\usepackage{multirow}
\usepackage{multicol}
\usepackage{algorithm}
\usepackage{appendix}
\usepackage{tikz}
\usepackage{pgfplots}
\newcommand{\mycomment}[1]{}
\usepackage[noend]{algpseudocode}

\title{Low-synchronization Arnoldi Methods for the Matrix Exponential with Application to Exponential Integrators}

\author{Tanya V. Tafolla, St\'{e}phane Gaudreault, Mayya Tokman}

\begin{document}

    \maketitle

    \section*{Abstract}
    High order exponential integrators require computing linear combination of exponential like $\varphi$-functions of large matrices $A$ times a vector $v$. Krylov projection methods are the most general and remain an efficient choice for computing the matrix-function-vector-product evaluation when the matrix is $A$ is large and unable to be explicitly stored, or when obtaining information about the spectrum is expensive. The Krylov approximation relies on the Gram-Schmidt (GS) orthogonalization procedure to produce the orthonormal basis $V_m$. In parallel, GS orthogonalization requires \textit{global synchronizations} for inner products and vector normalization in the orthogonalization process. Reducing the amount of global synchronizations is of paramount importance for the efficiency of a numerical algorithm in a massively parallel setting. We improve the parallel strong scaling properties of exponential integrators by addressing the underlying bottleneck in the linear algebra using low-synchronization GS methods. The resulting orthogonalization algorithms have an accuracy comparable to modified Gram-Schmidt yet are better suited for distributed architecture, as only one global communication is required per orthogonalization-step. We present geophysics-based numerical experiments and standard examples routinely used to test stiff time integrators, which validate that reducing global communication leads to better parallel scalability and reduced time-to-solution for exponential integrators.

    \section{Introduction}

        Many problems in science and engineering require simulating complex systems with dynamics governed by various forces acting over a wide range of time scales. Numerical solutions of such systems typically involve spatial discretization of the partial differential equations (PDEs) and their transformation into large systems of ordinary differential equations (ODEs) using the method of lines. Because these systems span wide temporal scales, the resulting ODEs are typically stiff, necessitating numerical solutions over intervals far exceeding the characteristic time scale of the fastest modes. Therefore, selecting an appropriate time integration method is crucial for ensuring both accuracy and efficiency in computing solutions.
        
        Explicit time integrators are unsuitable for such problems due to strict stability limitations on the timestep and implicit methods are often the preferred choice. With implicit schemes, stiffness manifests in the increased computational complexity of both nonlinear and linear solvers involved in approximating solution at each timestep.  A common numerical approach to implementing an implicit method is to use the Jacobian-Free Newton-Krylov algorithm \cite{knoll2004jacobian}. In this method the nonlinear system of equations resulting from the implicit discretization is handled with a Newton iteration which, in turn, involves  solution of linear systems using  Krylov subspace solvers. For stiff systems, the Krylov method can converge slowly, and a preconditioner is needed to make it efficient. However, constructing an effective preconditioner is often challenging and is the subject of considerable research.

        Over the past decades, exponential integrators have emerged as an efficient alternative for numerically integrating large-scale stiff systems. Exponential time integrators express the numerical solution as a linear combination of products of exponential-like functions of the Jacobian, or its approximations, and vectors \cite{johnthesis, valthesis, tokjcp, mayya_neweprk, hochbruck2010exponential}. Similar to implicit methods, exponential schemes have good stability properties that allow integration with a large timestep. For problems where a fast preconditioner is not available, exponential methods can offer computational savings per timestep compared to implicit integrators \cite{loffeld2013comparative,mayya_neweprk}. This is particularly true when Krylov-projection based algorithms are used for both implicit and exponential integrators. For many large-scale applications of interest, evaluating exponential-like functions of matrices can be more efficient and accurate than computing a rational function of a Jacobian, which is part of an implicit timestep \cite{loffeld2013comparative}.
        
        The dominating cost of exponential integrators is the computation of the exponential-like $\varphi$-functions times a vector. Krylov methods have shown to be efficient for problems involving large matrices that are difficult to compute explicitly or store, or when obtaining information about the spectrum is prohibitively expensive \cite{moler2003nineteen}. The main idea of Krylov subspace methods is that matrix functions can be efficiently computed using a projection of the matrix onto a low-dimensional subspace. For the general case where the matrix is non-symmetric, this is generally accomplished through the Arnoldi iteration. At the base of the Arnoldi method is essentially a Gram-Schmidt orthogonalization process which requires a projection of consecutive Krylov vectors onto an orthogonal basis computed during previous iterations. Thus, the  Arnoldi method does not scale well in parallel due to the cost of global communication needed to compute the inner products involved in the orthogonalization procedure. This global communication can become a bottleneck as the number of processors increases \cite{globalcomm_kimpact}.

        Various techniques have been proposed in recent work on the General Minimal Residual (GMRES) method  to reduce global communication by reformulating the orthogonalization process into a projection onto the orthogonal complement \cite{delayedcgs, wy, thomas2023iterated}. Such algorithms aim to decrease the amount of the global communication required by the Arnoldi iteration to a single call per Arnoldi-step, improving the strong scaling properties of linear solvers on distributed memory parallel systems. 
        
        The aim of this paper is to explore Krylov methods that minimize global communication and to extend their application to computing the matrix exponential, with a primary focus on enhancing the overall performance of exponential time integrators. We demonstrate the performance advantages of the new methods using a variety of large-scale problems, including some standard benchmarks routinely used by the numerical weather prediction community \cite{galewsky2004initial, swbenchmark}. Our numerical experiments indicate that these low-synchronization orthogonalization methods are strong competitors to current techniques when using a large number of processors. The methods also exhibit a clear trend of improved strong scaling behavior as the number of processors increases, resulting in enhanced efficiency for exponential time integrators. The paper explores the potential of these methods to address the communication bottleneck and improve the overall performance of exponential integrators for large-scale simulations. 
        
        The remainder of this paper is organized as follows: Section 2 provides background information and highlights parallel scaling limitations. Section 3 presents current low-synchronization variants and introduces a new hybrid low-synchronization Arnoldi algorithm. Section 4 shows numerical evidence of the improvement in parallel scalability achieved by different exponential time integrators when using low-synchronization Arnoldi methods. Finally, Section 5 draws conclusions and outlines future research directions.




    \section{Background and Motivation}




        Exponential time integrators are a family of numerical methods for solving initial value problems of the form 
			\begin{align}\label{eq:ode}
			\dfrac{d }{d t} u(t) = F(u(t)), \quad u(0) = u_0, \quad u(t) \in \mathbb{R}^N
			\end{align}
        \noindent where $u(t)$ is the solution of (\ref{eq:ode}) at time $t$, $N$ is the number of degrees of freedom, and $F(u(t))$ is the forcing term. If the time interval $[t_0, t_{end}]$ is discretized over a set of nodes $\{t_n \}_{n=0}^{M}$, then the solution of (\ref{eq:ode}) at time $t_n$ is given by $u(t_n)$. Using Taylor expansion about the exact solution $u(t_n)$, the system \eqref{eq:ode} can be written as
		 	\begin{equation}\label{eqn:ode_taylor_expand}
		 	\frac{d u}{d t} = F(u(t_n)) + J(u(t_n))(u - u(t_n)) + R(u(t)),
		 	\end{equation}
		 	
	\noindent  where $J(u(t_n)) = \dfrac{\partial F}{\partial u}(u(t_n))$ is the Jacobian matrix, $F(u(t_n))$ is the right-hand side of (\ref{eq:ode}) at time $t_n$, and $R(u(t)) = F(u(t)) - F(u(t_n)) - J(u(t_n))(u(t) - u(t_n))$ is the nonlinear remainder after two terms in the Taylor expansion. Employing $e^{-J_n t}$ as an integrating factor and integrating \eqref{eqn:ode_taylor_expand} exactly over the time interval $[t_n, t_{n+1}]$  yields the following Volterra integral equation

        \begin{equation}\label{eq:exp_int_form}
		 	u(t_{n+1}) = u(t_n) + \varphi_1(hJ(u(t_n))) h F(u(t_n)) + \int_{t_n}^{t_{n+1}} e^{(t_{n+1} - t)J_n}R(u(t)) dt,
		 	\end{equation}
    
        \noindent where $\varphi _{1}(z) = (e^{z} -1)/z$. Equation (\ref{eq:exp_int_form}) is  the starting point for deriving various exponential integrators by replacing the exact values with the numerical approximation, i.e. $u_n \approx u(t_n)$ or $J_n \approx \dfrac{\partial F}{\partial u}(u(t_n))$, and approximating the nonlinear integral in the right-hand-side to estimate $u_{n+1}$. For example, well-known quadratures that use either nodes $t_i$  or successively computed approximations of the solution over the interval $[t_n, t_{n+1}]$ will result in the multistep- or Runge-Kutta-type integrators respectively \cite{johnthesis,tokjcp, mayya_neweprk}. Any such polynomial approximation of the nonlinear remainder function $R(u(t))$ will result in the approximate solution $u_{n+1}$ being expressed as a linear combination of exponential-like $\varphi$-functions multiplied  by vectors
        \begin{equation}\label{eq:ansatz}
		\varphi_0(J_n)b_0 + \varphi_1(J_n)b_1 + \varphi_2(J_n)b_2 + ... + \varphi_p(J_n)b_p
        \end{equation}
        with the $\varphi$-functions defined as
		\begin{align}\label{eq:phi_int_def}
		\varphi_0(z) &= e^{z} \\
		\varphi_k(z) &= \int_0^1 \frac{e^{(1 - \theta)z}\theta^{k-1}}{(k-1)!}d\theta \quad \forall k \in \mathbb{N}^{> 0}.
		\end{align} 
     
        Evaluating linear combinations of the form \eqref{eq:ansatz} constitutes the main computational cost of exponential time integrators. Various approaches exist for computing the $\varphi$-functions, including Taylor series \cite{proofoneexp}, Leja methods \cite{caliari2016leja,deka2023lexint}, interpolation methods \cite{caliari2023bamphi}, and Krylov subspace methods \cite{kiops, koskela_iop, sidjeio}. We are interested in problems where matrix $J_n$ is large, stiff and might not even be available explicitly. For example, $F$ could be a result of a complex spatial discretization of the forcing terms in a stiff system of PDEs and an explicit analytic evaluation of its Jacobian matrix might be completely impractical from the efficiency and memory storage perspectives. Often only a function that computes a product of the Jacobian matrix with a vector is provided. The most appropriate numerical approach to evaluate linear combinations \eqref{eq:ansatz} in such applications becomes the Krylov subspace projection-based methods.  Improving parallel scalability of the Krylov projection is the focus of this paper.     

        \subsection{Linear combination of $\varphi$-functions} \label{sec:linearPhiFunc}            
            To optimize evaluation of linear combinations of $\varphi$-functions and vectors we follow previous work \cite{proofoneexp, kiops, sidje_expokit} and reformulate the linear combination \eqref{eq:ansatz} as a single matrix exponential using the equivalence
                
            \begin{equation}\label{eq:philincomb}
               \sum_{j = 0}^{p}\tau^j \varphi_j(\tau J_n)v_j = \begin{pmatrix} I_{N} & 0\end{pmatrix} \, \exp(\tau A) \, 
               \begin{pmatrix} v_0\\ e_p\end{pmatrix}, 
               \end{equation}
                    
                    \noindent where the augmented matrix $A$ is defined as 
        			\begin{equation}  
                    A = \begin{bmatrix}
        				 \tau J_n & B \\ 0 & K
        				  \end{bmatrix} \in \mathbb{R}^{(N+p)\times(N+p)},
                    \end{equation}

                    \begin{equation}
                        B = [b_p, ..., b_1] \in \mathbb{R}^{N\times p}, \, K = \begin{bmatrix}
        				 0 & I_{p-1} \\ 0 & 0 
        			    \end{bmatrix} \in \mathbb{R}^{p \times p}, \, \tau \in \mathbb{R},
                    \end{equation}
                    
                    \noindent and $e_p$ is the $p$th column of $I_p$, the identity matrix of size $p$.

                    In order to match the linear combination of the general ansatz \eqref{eq:ansatz}, the exponential of the augmented matrix must be computed with $\tau = 1$. It is important to note that assembling the augmented matrix $A$ is not required; instead, the implementation of this technique only necessitates the product of $J_n$ with a vector of compatible dimension. The task now is to consider how the exponential of the augmented matrix can be efficiently evaluated using a Krylov subspace method. 

                \subsection{Krylov Subspace Approximation}\label{sec:KryApprox}

                Krylov subspace approximations of the products of matrix functions and vectors of type $f(A)b$ are based on the projection of matrix $A \in \mathbb{R}^{d \times d}$ and nonzero vector $b \in \mathbb{R}^d$ onto a smaller Krylov subspace of dimension $m \ll d$,
                \begin{equation}
                \mathcal{K}_m(A,b) = \text{span} \{b,  Ab, A^2b, ... , A^{m-1}b \}.    
                \end{equation}
                In the most general case, when $A$ is a non-symmetric matrix, the Arnoldi process is employed to compute the Hessenberg reduction
                \begin{equation}\label{eq:Hdecomp}
                    H_m = V_m^TAV_m,
                \end{equation}
                where $V_m = [v_1, \ldots, v_m]$ is a matrix whose columns form an orthonormal basis of $\mathcal{K}_m$, and $H_m\in \mathbb{R}^{m \times m}$ is upper Hessenberg. Gram-Schmidt algorithm is commonly used as part of the Arnoldi process to generate a set of orthonormal vectors $ \{ v_1, v_2, ..., v_{m+1}\}$. Note that the matrix $V^T_m V_m$ is the identity when the basis  $ \{ v_1, v_2, ..., v_{m+1}\}$ is orthonormal. This property can be used to assess the departure from orthogonality in the numerical algorithms.  In matrix form Arnoldi iteration can be written as
                \begin{equation}\label{eq:arnoldidecomp1}
                AV_m = V_m H_m + h_{m+1,m}v_{m+1} e_m^T,
                \end{equation}
                where $e_m$ is the $m$th unit vector. Using Krylov projection the matrix exponential can be approximated as      	
                    \begin{equation}\label{eq:kryexp}
                        e^Ab \approx \lVert b \lVert_2 V_m e^{H_m}e_1.
                    \end{equation}
                The magnitude of the error associated with first evaluating the matrix exponential function on the smaller Krylov subspace $\mathcal{K}_m(A, b)$ and expanding the result back to the original space $\mathbb{R}^N$ is significant when $\lVert A \rVert_2$ is large (see, e.g., Theorem 4.3 in \cite{saaderror}). Thus, the convergence of the Krylov iteration is slow when $A$ is stiff. A common strategy to mitigate these issues is to use adaptive sub-stepping, analogous to the techniques used in many solvers for ordinary differential equations \cite{kiops,adaptkry, sidje_expokit,  sidjeio}. The basic idea, exemplified here for the computation of $w(\tau_{i+1}) = e^{\tau_iA}w(\tau_i)$, is to solve an equivalent linear initial value problem              
                \begin{equation}\label{eq:ode_substep}
                      \frac{d}{dt}w(\tau) = A w(\tau), \qquad w(0) = b. 
                  \end{equation}
            The action of the matrix exponential is obtained through the recurrent evaluation of the solution to this initial value problem, referred to as sub-stepping
                \begin{subequations}
                    \begin{alignat}{2}
                      e^{\tau A}b &= e^{(\tau_n + ... + \tau_2 + \tau_1 )A}b \\
                    &= e^{\tau_n A}(...(e^{\tau_2 A}(e^{\tau_1 A}b)),
                       \label{subeqn:substep}
                    \end{alignat}
                  \end{subequations}
                \noindent where each product $e^{\tau_i A}b_{\tau_i}$ in (\ref{subeqn:substep}) is approximated with a Krylov projection (\ref{eq:kryexp}).  
                
                This sub-stepping procedure reduces the norm of the matrix and speeds up convergence of the Krylov iteration. At each sub-step, both the size $\tau_i$ and Krylov subspace dimension $m_i$ are selected adaptively as in \cite{kiops}, and a maximum threshold on $m_i$ is set based on the characteristics of the hardware.

                \subsection{Parallel scaling limitations of the Krylov subspace projection}

                As outlined above, the sub-stepping procedure requires an Arnoldi iteration to be performed at each substep $\tau_i$.  The overall computational cost of an exponential time integrator that utilizes the Krylov subspace projections to evaluate linear combinations \eqref{eq:ansatz} is then primarily dependent on the efficiency of the underlying Arnoldi iterations. Using this algorithm on a distributed memory computational platform usually assumes that vector entries are distributed among $P$ processors. Some computations such as vector scaling or vector addition can be done locally by each individual processor. However, other computations such as inner products require global communication. Each processor handles its respective portion, then the intermediate results are shared among the remaining processors, added through global communication, and broadcasted back to each processor, for instance in the form of a \texttt{MPI\_Allreduce} call\footnote{The MPI\_Allreduce function of the Message Passing Interface (\texttt{MPI}) is a collective communication operation that performs a global reduction operation (for example a summation) on a set of input values and distributes the result to all participating processes. This function is blocking and introduces collective synchronization into the exponential solver.}. 
                
                The primary challenge in implementing the Arnoldi process for parallel computational platforms lies in the Gram-Schmidt orthogonalization. Two common methods for orthogonalization are Modified Gram-Schmidt (MGS) and Classical Gram-Schmidt (CGS). Although CGS requires only two global communication calls at each $j^{th}$ Arnoldi step (one for all the inner products and one for vector normalization), this algorithm is less stable and more susceptible to rounding errors compared to MGS. MGS is a more stable algorithm, but is significantly more expensive in parallel, requiring  $(j+1)$ global communications ($j$ inner products and one for vector normalization, see Algorithm \ref{al:mgs} ) at the $j^{th}$ Arnoldi step. For efficient Krylov-projection approximation of $e^{\tau_i A}b_i$ it is therefore desirable to develop an algorithm that reduces the total global communication during the orthogonalization process without compromising the stability and the accuracy of the approximation.
                
                \begin{algorithm}[H]
			\caption{Arnoldi proces (with modified Gram-Schmidt)}
			\begin{algorithmic}[1]
				\Procedure{Arnoldi}{$A,b, m$}       \Comment{matrix, initial vector, size of subspace}
				\State $\beta = \lVert b \lVert_2$
				\State $v_1 = b/ \beta$
				\State $j = 0$
				\While{$j < m$}
				\State j = j+1
				\State $v_{j+1} = Av_{j}$ \label{comm1}\Comment{next vector in basis}
                    \For{$i = 1...j$} 
                    \Comment{orthogonalization}
                    \State $H_{ij} = v_i^Tv_{j+1}$ \Comment{global sync}\label{comm2}
                    \State $v_{j+1} = v_{j+1} - H_{ij}v_i$
                    \EndFor
                    \State nrm = $\lVert v_{j+1} \lVert_2$ \Comment{global sync}\label{comm3}
				\If{nrm $\approx$ 0} \Comment{check for lucky breakdown\cite{luckybreakdowns}}
				\State break
				\EndIf
				\State $v_{j+1} = v_{j+1}/\texttt{nrm}$  
				\State $H_{j+1,j} = \texttt{nrm}$
				\EndWhile			
				\EndProcedure
			\end{algorithmic}\label{al:mgs}
		\end{algorithm}

           \subsubsection{Incomplete orthogonalization process}\label{sec:io}
                One way to reduce the total amount of global communication is the Incomplete Orthogonalization Process (IOP), in which new vectors are orthogonalized relative to a few previous basis vectors, e.g., the previous two,  as opposed to all of them. Incomplete orthogonalization was previously used for solving linear systems and eigenvalue problems \cite{saad_og_io_eig, saad_og_io} and has since been extended to matrix exponential computations \cite{kiops,koskela_iop, sidjeio}. 
                
                While the IOP seems to improve overall parallel scalability and works well in many practical applications (e.g. \cite{val_cube}), its behavior is difficult to analyze. In contrast with the standard Arnoldi procedure, the IOP produces a $H_m$ which is not Hessenberg, the columns of $V_m$ do not form an orthonormal basis and $V^T_m\,V_m$ is not the identity matrix. Furthermore, in theory, rounding errors can result in basis vectors that are nearly collinear with those not orthogonalized, potentially compromising the accuracy of the approximation \eqref{eq:kryexp}. In some applications IOP does not produce sufficient accuracy and full orthogonalization is required \cite{stewart2023variable}.  Thus, it is desirable to develop a Krylov method that uses complete orthogonalization that can be efficiently scaled for parallel computing platforms.

                For our study, we use the \texttt{KIOPS} algorithm, a Krylov Incomplete Orthogonalization Procedure developed by Gaudreault et al. that combines incomplete orthogonalization (orthogonalizing against the previous two vectors) from \cite{koskela_iop, sidjeio}, substepping of $e^Av$, and Krylov adaptivity in Section \ref{sec:KryApprox} of the augmented matrix from Section \ref{sec:linearPhiFunc}. This algorithm has substepping and adaptivity techniques for $\tau_i$ and $m_i$ that do not require multiplication with the matrix $A$ and take into account the impact of round-off error, resulting in an overall highly efficient and accurate algorithm for computing $e^{\tau_i A}b_i$ \cite{kiops}. This method was shown to improve the efficiency of exponential integrators for solving commonly tested stiff PDEs and applications ranging from combustion to weather prediction \cite{val_cube, kiops, stewart2023variable}. The structure of the algorithm (further discussed in Section \ref{sec:expint}) also allows for easy comparisons of the orthogonalization schemes that will be presented in the following section.

        \section{Low-synchronization Arnoldi Algorithm for Krylov Approximations}               






        Arnoldi algorithms are used for a wide variety of applications ranging from finding eigenvalues (and corresponding eigenvectors) \cite{delayedcgs,hernandez, kimchronolagnorm} to solving linear systems $Ax = b$ with GMRES \cite{wy, hidelatency, delayedcgs, postmodern, sstep_thomas}. The parallel efficiency of these algorithms has been a subject of extensive interest and research, dating back to the mid 90's \cite{sidje1997alternatives,philippe1995parallel,de1994communication}. The most recent developments in improving the parallel scalability of Krylov linear solvers on massively parallel computing platforms have been focused on developing low-synchronization Arnoldi methods \cite{delayedcgs,hidelatency, wy, postmodern, sstep_thomas}. Low-synchronization methods aim to reduce the amount of global communication or synchronization points in an algorithm, such as \texttt{MPI\_Allreduce} calls. Low-synchronization Arnoldi methods are derived by reformulating the orthogonalization process of Arnoldi algorithm to trade off increased computational complexity for reduced global communication needed at each Arnoldi step. In this study, we focus on Modified Gram-Schmidt variants of the low-synchronization methods that require one global synchronization per Arnoldi-step. We consider the \texttt{WY} methods from \cite{wy, ThomasNeumann} and the iterated Gau{\ss}-Seidel formulations \cite{postmodern}. These algorithms employ reformulation of the orthogonalization process in MGS by expressing it as the projection onto the orthogonal complement. The following sections describe the low-synchronization methods and their extension to exponential integrators.   

        \subsection{Main ideas of low-synchronization methods}

        \subsubsection{Projection onto the orthogonal complement}\label{sec:introprojortho}

        A key idea of many modern low-synchronization techniques is that the orthogonalization of each subsequent vector $v_{j+1}$ in the Krylov basis at Arnoldi step $j$ is done by projecting it onto the orthogonal complement of the $j$-dimensional Krylov subspace $V_j^{\perp} = \{ w \in \mathrm{R}^n : \langle v, w \rangle = 0 \text{, } \forall v \in V_j \}$ \cite{ruheogwy, wy} . The projection matrix onto $V_j^{\perp}$ is given by 
        \begin{equation}\label{eqn:orthoproject} 
            P = I_j - V_j(V_j^TV_j)^{-1}V_j^T
            \end{equation} and $v_{j+1}$ can be orthogonalized as follows
            \begin{equation}\label{eqn:projectorthocomp}
			 Pv_{j+1} = (I_j - V_j(V_j^TV_j)^{-1}V_j^T)v_{j+1}
	\end{equation}
        where $I_j$ is the $j \times j$ identity matrix. The sequential orthogonalization of the MGS algorithm (lines 8 - 10 of Algorithm \ref{al:mgs}) can then be replaced by a matrix-vector multiplication with the projection matrix $P$ and vector $v_{j+1}$. The interior matrix $(V_j^TV_j)^{-1}$ can be used to correct for the non-orthogonality of the basis due to round-off errors in finite precision arithmetic; although the inner products $v_i^Tv_j$ are expected to be small, they are not exactly zero. 
		
        The low-synchronization Arnoldi algorithms explored in this study are rooted in similar principles but differ in terms of possible loss of orthogonality and computational complexity that vary based on the distinct approximations to the correction matrix $T_j \approx (V_j^TV_j)^{-1}$. The different approaches to estimating the correction matrix $T_j$ are discussed in the subsequent sections.
           
        \subsubsection{Lagged normalization}\label{sec:stablagnrm}
           
        The second strategy used to minimize communication is lagged normalization \cite{hernandez,kimchronolagnorm}. The key idea is to postpone the normalization of the next vector $v_{j+1}$ until the communication needed for the norm can be grouped with other calculations such as inner products in the subsequent step. In this paper, we use low-synchronization methods each with their unique approximation to the orthogonal projector along with lagged normalization. 

        \subsection{Hybrid low-synchronization}






        The lagged normalization strategy delays the normalization of vector $v_{j+1}$ at Arnoldi step $j$. Grouping the computation of the norm for $v_{j+1}$ at the next Arnoldi step with other inner products reduces the total amount of global communication to one per-Arnoldi-iteration, i.e. eliminating the global communication from the vector normalization. Although lagged normalization is essential for the orthogonal complement-based orthogonalization to remain a one-synchronization method, the vector $v_{j+1}$ will potentially have a large norm. This will result in an unstable iteration when the matrix operator is applied to $v_{j+1}$.

        An alternative would be to introduce a Pythagorean form of a norm estimate \cite{hidelatency,hernandez}. At the end of each Arnoldi step $j$, the norm estimate of $v_{j+1}$ can be computed and $v_{j+1}$ can then be ``normalized" using the following estimate \cite{hidelatency,hernandez}:

            \begin{equation}\label{eq:nrmest}
                    \texttt{nrm} = \sqrt{||v_{j+1}||_2^2 - \sum_{i=0}^j (v_i^Tv_{j+1})^2 }.
                \end{equation}
                The matrix operator applied to $v_{j+1}$ at the subsequent iteration should then be stable. Furthermore, the inner products $V_{1:j+1}^Tv_{j+1}$ required by the norm estimate \eqref{eq:nrmest} are obtained by the grouped communication from section \ref{sec:introprojortho} (details in sec. \ref{sec:cwy}). Therefore, no additional communication is required, and the algorithm remains a one-synchronization method. However, this norm estimate is not exact for MGS and can be sensitive to round-off error if the terms inside the square root are similar in magnitude, i.e if orthogonality of $V_m$ is lost \cite{hidelatency}. This can cause the argument of the square root to become negative and the algorithm to break down. The inaccuracy of the norm estimate can also impact the accuracy of the Krylov projection since the norm values are entries of the Hessenberg matrix used in computing the projection of the matrix-exponential (eq. \ref{eq:kryexp}). Therefore, a combination of the orthogonal complement projection and norm estimate is also potentially unstable.

        Given the instabilities of the orthogonal complement-based orthogonalization paired with either lagged normalization or the norm estimate, we present a hybrid version of low-synchronization algorithms that combines all three techniques (projection to the orthogonal complement, lagged normalization and norm estimate); an outline of the hybrid low-synchronization algorithm is presented in Algorithm \ref{alg:generalhybrid}. At the end of Arnoldi step $j$, the norm estimate of vector $v_{j+1}$ is computed using equation (\ref{eq:nrmest}), allowing for a stable matrix-vector product at the following step. At the next step, after the matrix-vector product to extend the Krylov basis, the vectors $v_{j+1}$ and $v_j$ (current and previous vector) are re-scaled by the norm-estimate of $v_j$, which was computed at the previous iteration ($V_{j:j+1} = \texttt{nrm\_est}*V_{j:j+1}$). The re-scaling of $v_j$ will ensure that the true norm is obtained from the group computation of inner products (details in sec. \ref{sec:cwy}). The true norm of the previous vector is then used in the Hessenberg matrix, and to accurately normalize the previous vector. This new hybrid low-synchronization method requires only one global synchronization, avoids multiplying matrix $A$ and vector $v_j$ with a large norm and uses the true norm obtained from lagged normalization, avoiding inaccuracies from finite precision stemming from the norm estimate. We call these methods ``hybrid' because we incorporate lagged normalization and perform normalization at the end of each Arnoldi step. 

        Sections \ref{sec:cwy} to \ref{sec:gsmgs} expand on details regarding the grouped computation for obtaining all necessary terms (line 6, Algorithm \ref{alg:generalhybrid}), each low-synchronization method's approximation to the correction matrix $T$ (line 7, Algorithm \ref{alg:generalhybrid}), and the inclusion of the norm-estimate computation (line 5 and 9, Algorithm \ref{alg:generalhybrid}) .
            
        \begin{algorithm}[H]
			\caption{Outline of hybrid low-synchronization algorithm with techniques (sec. \ref{sec:introprojortho}), (sec. \ref{sec:stablagnrm}), and norm estimate (eq. \ref{eq:nrmest})}
			\begin{algorithmic}[1]
				\Procedure{hybrid\_lowsync}{$A, v_1, m$}       \Comment{matrix operator, initial vector, size of subspace}
				\For {$j = 1,...m$} 
				\State Expand Krylov basis $\{v_1, v_2,...\}$ : $v_{j+1} = Av_{j}$
                    \If{$j > 1$}
                    \State Re-scale $v_{j+1}$ and $v_j$ by norm estimate: $v_{j:j+1} = \texttt{nrm\_est} \times v_{j:j+1}$
                    \EndIf
				\State Grouped computation for projection (\ref{eqn:projectorthocomp}); includes lagg normalization and norm est (\ref{eq:nrmest}).
                    \State Formation of correction matrix $T \approx (V^TV)^{-1}$%
				\State Projection via (\ref{eqn:projectorthocomp})
                    \State ``Normalization'' of $v_{j+1}$ via (\ref{eq:nrmest})
				\EndFor			
				\EndProcedure
			\end{algorithmic}\label{alg:generalhybrid}
	\end{algorithm}
           
	\subsubsection{Hybrid Compact WY MGS}\label{sec:cwy}

        Recall the orthogonalization for-loop in MGS (Algorithm \ref{al:mgs_orthoonly}), here the superscript indexing $i$ is added to explicitly label $v_{j+1}$ at the $i^{th}$ step of the recursive orthogonalization.
        \begin{algorithm}[H]
			\caption{Modified Gram-Schmidt Orthogonalization}
			\begin{algorithmic}[1]
				\Procedure{mgs\_ortho}{$A,b, m$}       
                    \For{$i = 1...j$} 
                    \Comment{orthogonalization}
                    \State $H_{ij} = v_i^Tv_{j+1}^{(i)}$ \Comment{global sync}\label{comm2}
                    \State $v_{j+1}^{(i+1)} = v_{j+1}^{(i)} - H_{ij}v_i$
                    \EndFor
				\EndProcedure	\end{algorithmic}\label{al:mgs_orthoonly}
		\end{algorithm} 
        \noindent If the for-loop is unrolled at iteration $i=j$, we can obtain an expression for $v_{j+1}^{(i+1)}$ in terms of $v_{j+1}^{(1)} = Av_{j}$. As an illustrative example, consider the for-loop expansion for $j=2$ (orthogonalizing $v_3$ against $v_1$ and $v_2$)

        \begin{subequations}
            \begin{alignat}{2}
                v_{3}^{(3)} &= v_{3}^{(2)} - \langle v_2, v_{3}^{(2)}\rangle v_2 \\
            v_3^{(3)} &= v_3^{(1)} - \langle v_1, v_3^{(1)} \rangle v_1 - \langle v_2, v_3^{(1)} - \langle v_1, v_3^{(1)} \rangle v_1 \rangle v_2 \label{eqn:unroll_mgsloop}
            \end{alignat}
        \end{subequations} The expansion (\ref{eqn:unroll_mgsloop}) can be simplified and re-written as the application of the rank-1 projectors $(I - v_jv_j^T)...(I-v_1v_1^T)v_{j+1}$\cite{hernandez,wy}. The product of these rank 1 projectors can then be expressed in matrix form as 		\begin{equation}\label{eqn:wyproj}
			(I - v_jv_j^T)...(I-v_1v_1^T)= I - V_j T_j V_j^T.
			\end{equation}
        The matrix form of the rank-1 projections (\ref{eqn:wyproj}) then forms an approximation to the orthogonal complement, with the lower triangular loss-of-orthogonality (or correction) matrix $T_j$ approximating the interior symmetric matrix of the projection $(V^TV)^{-1}$ from (\ref{eqn:orthoproject}). The MGS-Arnoldi process is then re-written as a projection to the orthogonal complement through the application of (\ref{eqn:wyproj}) instead of the sequential for-loop from Algorithm \ref{al:mgs_orthoonly}. This re-formulation of the MGS-Arnoldi through the orthogonal complement is based on Ruhe \cite{ruheogwy}, and has been recently extended by \'{S}wirydowicz et al. employing lagged normalization; we will refer to it as the \textit{Compact WY} MGS (Algorithm 6 in \cite{wy}). The correction matrix $T_j$ is computed one row at a time at each Arnoldi step $j$ through the recursion formula 

            \begin{subequations}\label{eqn:substep}
                    \begin{alignat}{2}
                      T_1 &= [1] \\
			T_j &= \begin{bmatrix}
			T_{j-1} & 0 \\
			 -(V_{j-1}^Tv_j)^T T_{j-1}& 1
			\end{bmatrix}. \label{eqn:recursionForT}
                    \end{alignat}
            \end{subequations} Through this re-formulation only one global reduction sum is needed at each Arnoldi step for the computation of $V_{1:j}^T V_{j:j+1}$ \cite{wy}. The resulting $(j \times 2)$ matrix is given by 
            \begin{equation}\label{eq:output_wy}
            \begin{bmatrix}
            v_1^T v_j & v_1^T v_{j+1} \\
            v_2^T v_j & v_2^T v_{j+1} \\
            \vdots & \vdots \\
            v_j^T v_j & v_j^T v_{j+1} \\
            \end{bmatrix}.
            \end{equation} After computation (\ref{eq:output_wy}), we have all components necessary for the projection.  The second column of (\ref{eq:output_wy}) is $V^T_jv_{j+1}$ from the projection  (\ref{eqn:projectorthocomp}), the first $(j-1)$ rows of column one are used for in the recursion formula for $T_j$ (\ref{eqn:recursionForT}), and the square norm of the previous vector $v_j$ is row $j$ of column one. 
            
            However, in order to incorporate the norm estimate, $||v_{j+1}||_2^2$ is needed. Therefore, we extend the computation (\ref{eq:output_wy}) to $V_{1:j+1}^TV_{j:j+1}$. This will result in a $(j+1) \times 2$ matrix whose entries are given by
            \begin{equation}\label{eq:output_wywnrm}
            \begin{bmatrix}
            v_1^T v_j & v_1^T v_{j+1} \\
            v_2^T v_j & v_2^T v_{j+1} \\
            \vdots & \vdots \\
            v_j^T v_j & v_j^T v_{j+1} \\
            v_{j+1}^T v_j & v_{j+1}^T v_{j+1}
            \end{bmatrix}.
            \end{equation} with the additional component of $||v_{j+1}||_2^2$ (row $(j+1)$ of the second column) for the norm estimate.  
            
            Once (\ref{eq:output_wywnrm}) is calculated, the remaining steps for the orthogonalization are computing the local matrix-vector product in the recursion formula for $T_j$ (\ref{eqn:recursionForT}) and the projection of $v_{j+1}$ via equation (\ref{eqn:projectorthocomp}). The values of the upper Hessenberg matrix $H_{1:j,j}$ are obtained by the local matrix-vector product $T_j V_j^Tv_{j+1}$. At the end of the Arnoldi step, $v_{j+1}$ is ``normalized" using equation (\ref{eq:nrmest}), and at the subsequent step, the computation $Av_j$ should be stable. In order to then obtain the true norm of the previous vector from computation (\ref{eq:output_wywnrm}) (i.e. row $j$, column $1$), vectors $V_{j:j+1}$ are re-scaled by the norm estimate (computed at the previous iteration) after the matrix-vector product $Av_{j}$ is computed (see line 15 in Algorithm \ref{alg:lowsync_skeleton}).
                                   
            We refer to the combination of the one-synchronization version (projection with lagged normalization) with a norm estimate computation as the hybrid \texttt{CWY} algorithm (\texttt{H-CWY}). The algorithm for a general hybrid low-synchronization method is presented in Algorithm \ref{alg:lowsync_skeleton} with the corresponding technique for obtaining the correction matrix in Algorithm \ref{alg:Tcwy}. The hybrid version remains a one-synchronization algorithm (line 13 of Algorithm \ref{alg:lowsync_skeleton} for the inner products of the grouped communication) but at the cost of additional computation from the norm estimate and a larger communication for the additional row of $V^TV$. 
                       
			\subsubsection{Hybrid Neumann Series MGS} \label{sec:neumannseries}

                The \texttt{CWY} method requires two matrix-vector products with computational complexity of $\mathcal{O}(m^2)$, which results in a higher computational cost than \texttt{MGS} \cite{wy}. In order to reduce this computational cost, the alternative method \textit{Truncated Neumann Series} MGS (\texttt{NWY}) presented by Thomas et al. (algorithm 4.1 in \cite{ThomasNeumann}) approximates the correction matrix $T$ with a finite Neumann series
                \begin{equation}\label{eqn:Nseries_pt1}
                    T_m = (I + L_m)^{-1} = I - L_m + L_m^2 - ... + L_m^p
                \end{equation}
                where $L_m$ is a strictly lower triangular matrix whose entries are
                \begin{equation}
                    L_{ij} =\begin{cases}
                        v_i^Tv_j, &\text{if } i> j \\
                        0 & \text{if } i\leq j 
                    \end{cases}.
                \end{equation}

                The Neumann series can be used to express the inverse of $(I-A)$ if  $\lim_{k\rightarrow \infty} A^k = 0$ (or equivalently, if the spectral radius $\rho(A) < 1$) as \cite{StewartNSeries}
			\[ (I-A)^{-1} = \sum_{k = 0}^\infty A^k.\]
                The inverse of a general matrix $Y$ can then be found by defining $A = I - Y$
                \begin{subequations}\label{eqn:epi6}
                \begin{alignat}{2}
                (I-A)^{-1} &= \sum_{k = 0}^\infty A^k \\
                    Y^{-1} &= \sum_{k = 0}^\infty (I - Y)^k.
                \label{subeqn:neumann_general_inverse}
                \end{alignat}
              \end{subequations}

                For the \texttt{CWY} method the inverse of $T_m$ is given by $T_m^{-1} = I + L_m$. Therefore, to find the correction matrix $T_m$, a Neumann series is used to approximate the inverse of $T_m^{-1}$ (i.e. $T_m$) using (\ref{subeqn:neumann_general_inverse}), with $A = (I - T_m^{-1}) = (I - (I + L_m)) = -L_m$.
                
                As an illustrative example, $T_3$ is shown below for the third iteration ($j = 3$) of the Arnoldi method, where $v_4$ is orthogonalized against $\{ v_1, v_2, v_3\}$
			
			\[ T_3 = \begin{bmatrix}
				1 & 0 & 0 \\
				-\alpha_{21} & 1 & 0 \\
				-\alpha_{31} + \alpha_{32}\alpha_{21} & -\alpha_{32} & 1
			\end{bmatrix}\] and $\alpha_{ij} = v_i^Tv_j$. In general the values of $T_m$ contain inner products of the columns of matrix $V_m$ and higher order terms containing the product of inner products ( i.e. $\alpha_{32}\alpha_{21}$ ). An approximation of $T_m$ can then be defined through a truncated Neumann series for $T_m^{-1}$:

                \begin{subequations}\label{eqn:epi6}
                \begin{alignat}{5}
                    (I - A)^{-1} &= \sum_{k = 0}^\infty A^k \\
			(I - (I - T_m^{-1}))^{-1} &= \sum_{k = 0}^\infty (I - T_m^{-1})^k \\
			T_m &= \sum_{k = 0}^\infty (I_m - T_m^{-1})^k = -L_m^k  \\
			T_m &\approx \sum_{k = 0}^2 -L_m^k \\
			T_m &\approx I_m - L_m^1 + L_m^2 \label{suneqn:nuemann_3terms}
                \end{alignat}
              \end{subequations}

            Following (\ref{suneqn:nuemann_3terms}) $T_3$ can then be written as 
            \begin{subequations}\label{eqn:epi6}
                \begin{alignat}{2}
                    T_3 &= \underbrace{ \begin{bmatrix}
				1 & 0 & 0 \\
				-\alpha_{21} & 1 & 0 \\
				-\alpha_{31}  & -\alpha_{32} & 1
				\end{bmatrix}}_{\approx T_3} +  \begin{bmatrix}
			0 & 0 & 0 \\
			0 & 0 & 0 \\
			\alpha_{32}\alpha_{21} & 0 & 0
			\end{bmatrix} \\
                T_m &\approx I - L_m.\label{subeqn:neumann_approx_t}
                \end{alignat}
              \end{subequations}
			
   
            If one defines $L^0$ to be the identity matrix, then the correction matrix $T_3$ can be expressed with three terms in the series expansion. The full correction matrix $T_m$ can be expressed as the sum of a lower triangular matrix with $\alpha_{ij}$ terms and matrices with the higher-order terms. Truncating the series with two terms will lead to a computationally cheaper approximation of $(V_j^TV_j)^{-1}$ that only includes the inner products $v_i^Tv_j$; the local matrix-vector product within the recursion formula (\ref{eqn:recursionForT}) is then avoided. 
            
            This method requires one global synchronization per Arnoldi-iteration for the same computation in (\ref{eq:output_wywnrm}) as the \texttt{CWY} method, and the hybrid version extends naturally. The only difference is the formation of the correction matrix $T_j$, which now avoids the additional matrix-vector product. The first $(j-1)$ rows of the first column of (\ref{eq:output_wywnrm}) will be the $j^{th}$ row of $L_{ij}$. Once $T_j$ is obtained, the vector $v_{j+1}$ can be orthogonalized by projecting onto the orthogonal complement via (\ref{eqn:projectorthocomp}) with the approximation of the correction matrix $T$ defined in (\ref{subeqn:neumann_approx_t}). The norm estimate of $v_{j+1}$ is computed at the end of the $j^{th}$ step, but the true norm is obtained at the following step through computation (\ref{eq:output_wywnrm}) (row $j$ of column one) after $V_{j:j+1}$ are re-scaled by the norm estimate. The values of the upper Hessenberg matrix $H_{1:j,j}$ are given by the local matrix-vector product $T_jV_j^Tv_{j+1}$ from the projection. 
            
            The algorithm for \texttt{H-NCWY} follows the outline from Algorithm \ref{alg:lowsync_skeleton} with Algorithm \ref{alg:Tncwy} for computing the correction matrix.
                  			
			\subsubsection{Hybrid Iterated Gau{\ss}-Seidel MGS} \label{sec:gsmgs}	

                
                The correction matrix $T_j$ for \texttt{CWY} and \texttt{NCWY} approximates the interior symmetric matrix $(V_m^TV_m)^{-1}$ with a lower triangular matrix. An alternative version to the \texttt{WY} methods previously presented is to work directly with the formula for the projection matrix and apply a two-step Gau{\ss}-Seidel iterative method to the normal equations within. The resulting $T_j^{(2)}$ correction matrix with two-steps will be almost symmetric and a better approximation of $(V^TV)^{-1}$. The method presented here is a simplified version of the technique originally developed for GMRES by Thomas et al. \cite{thomas2023iterated}, that requires only one global synchronization.
   
			The projection onto the orthogonal complement at Arnoldi step $j$ is given by
			\begin{equation}
			         P v_{j+1} = v_{j+1} - V_jx ,
			\end{equation}
   where the vector $x$ is the solution of the normal equations within the projection onto the orthogonal complement (\ref{eqn:normal_eqn_proj})
			   \begin{subequations}\label{eqn:gs1}
                \begin{alignat}{2}
				P v_{j+1} &= v_{j+1} - V_{j}\underbrace{(V_j^TV_j)^{-1}V_{j}^Tv_{j+1}}_{x} \label{eqn:normal_eqn_proj} \\
				V_j^TV_jx &= V_j^Tv_{j+1}. \label{eqn:projector_normaleqn}
				\end{alignat}
			\end{subequations} To solve (\ref{eqn:projector_normaleqn}), a two-step Gau{\ss}-Seidel relaxation scheme is used with the matrix splitting $ V^T_jV_j = M_j - N_j$, where $M_j = I + L_j$, $N_j = -L_j^T$ and $L_j$ is defined in section \ref{sec:neumannseries}. The choice of matrix splitting $M_j$ and $N_j$ was established by Ruhe when it was discovered that MGS employs a multiplicative Gau{\ss}-Seidel relaxation scheme \cite{ruheogwy} (see details in Appendix A of \cite{wy}), and further studied in Thomas et. al in \cite{postmodern} with the application to GMRES.
   
            The derivation of the 2-step Gau{\ss}-Seidel iteration to the normal equations is as follows:
			\begin{subequations}\label{eqn:gs2}
				\begin{alignat}{5}
				V_j^T V_j x &= \underbrace{V_j^T            v_{j+1}}_{b_j}\\
				(M_j - N_j)x &= b_j \\
				M_j x^{k+1} &= b_j + N_j x^{k}\\
				x^{k+1} &= M_j^{-1}(b_j + N_j x^{k}) \label{eqn:gs_iterative}
				\end{alignat}
				
			\end{subequations}
			
			Applying two steps of (\ref{eqn:gs_iterative}) with initial guess $x^0 = \textbf{0}$ yields the following correction matrix $T_j$:
			\begin{align*}
			x^1 &= M^{-1}b_j\\
			x^2 &= M^{-1}(b_j + N_j x^1)\\
			x^2 &= \underbrace{M^{-1}[I + NM^{-1}]}_{T_j}b_j
			\end{align*}
			Once the correction matrix $T_j$ is obtained, next vector in the subspace is then projected to the orthogonal complement, similarly to the WY methods, via equation (\ref{eqn:projectorthocomp}).
		
            This method requires one global synchronization, for the same computation $V_{1:j+1}^TV_{j:j+1}$ as in (\ref{eq:output_wywnrm}), and the hybrid version extends naturally. The first $j$ rows of the second column are the terms $V_{j}^Tv_{j+1}$ that are needed for the projection (\ref{eqn:projectorthocomp}). The $j^{th}$ element of the first column of (\ref{eq:output_wywnrm}) is the true square norm of the previous vector (obtained after the re-scaling of $V_{j:j+1}$, as in the \texttt{WY} methods). The first $(j-1)$ rows of the first columns are used for the matrices $N_j$, $M_j$, and $M_j^{-1}$; both $M_j$ and $N_j$ rely on matrix $L_j$, which are the first $(j-1)$ rows of $V_{1:j+1}^TV_{j:j+1}$ and $M_j^{-1}$ is the correction matrix from \texttt{CWY} in section \ref{sec:cwy}. The correction matrix $T_j$ is applied in two steps, first the matrix-vector product $a = [I + NM^{-1}]b_j$, followed by a lower triangular solve $Mx = a$. The values for the upper Hessenberg matrix $H_{1:j,j}$ are given by $T_jV^T_jv_{j+1}$. At the end of the orthogonalization, the norm estimate for $v_{j+1}$ is computed using (\ref{eq:nrmest}), and at the following step, the norm estimate is used to re-scale $V_{j:j+1}$ in order to obtain the true norm from computation (\ref{eq:output_wywnrm}).  The hybrid one-synchronization \textit{Iterated Gau{\ss}-Seidel MGS} (\texttt{H-GSMGS}) is outlined in Algorithm \ref{alg:lowsync_skeleton} with Algorithm \ref{alg:Tgsmsg} for computing the correction matrix. \newline


        \subsection{Low-synchronization algorithms for the matrix exponential}

            The low-synchronization methods presented above are based on the same principles in terms of the reformulation of the orthogonalization process paired with a lagged normalization to group inner products in a way that reduces the global communication to one per Arnoldi step. We extend the low-synchronization methods to include a norm estimate in order to avoid multiplying matrix $A$ with a vector $v_j$ with a large norm. We then compute the true norm and normalize at the subsequent step (referred to as hybrid). The key difference between each algorithm presented is the choice of the approximation to the interior symmetric matrix $(V^TV)^{-1}$. An outline of a hybrid low-synchronization method is presented in Algorithm \ref{alg:lowsync_skeleton} where the details are provided for the global communication, and additional scaling and computation needed for the hybrid algorithms. At line 24, this algorithm can be paired with Algorithms \ref{alg:Tcwy}, \ref{alg:Tncwy}, or \ref{alg:Tgsmsg} for computing a correction matrix $T$ and obtaining the values for the upper Hessenberg matrix $H_m$. 

            Once the components $V_m$ and $H_m$ of the Krylov subspace are obtained from Algorithm \ref{alg:lowsync_skeleton}, the matrix-exponential times a vector can be computed using the Krylov projection formula (\ref{eq:kryexp}). If $A$ (line 12 of Algorithm \ref{alg:lowsync_skeleton}) is the augmented matrix from section \ref{sec:linearPhiFunc} then one matrix function evaluation is needed within an exponential integrator to obtain all $\varphi$-functions. 


            \begin{algorithm}[H]
			\caption{Hybrid low-synchronization Arnoldi }
			\begin{algorithmic}[1]
				\Procedure{HGSMGS}{$m,b$}       \Comment{size of subspace, initial vector}
                \If{\texttt{H-CWY} or \texttt{H-NCWY}}
                    \State $T_{1:m,1:m} = I_m$ \Comment{initialize correction matrix to identity}
                \Else \Comment{if using \texttt{H-GSMGS}}
				\State $M_{1:m,1:m} = I_m$ \Comment{initialize to identity matrix}
				\State $M^{-1}_{1:m,1:m} = I_m$ \Comment{initialize to identity matrix}
				\State $N_{1:m,1:m} = 0_m$ \Comment{initialize to zero matrix}
                \EndIf
				\State $v_1 = b$
				\State $j = 0$
				\While{$j < m$}
				\State j = j+1
				\State $v_{j+1} = Av_{j}$ \Comment{next vector in basis}
				\State $\texttt{temp} = V^T_{1:j+1}V_{j:j+1}$ \Comment{\texttt{MPI\_AllReduce}; computation (\ref{eq:output_wywnrm})}
                    \If{$j > 1$}
                    \State $V_{j:j+1} = V_{j:j+1}*\texttt{nrm\_est}$ \Comment{re-scale by previous norm estimate}
                    \EndIf
                    \State $\texttt{true\_nrm} = \sqrt{\texttt{temp}_{j,1}}$ \Comment{true norm of prev. vector}
                    \If{$j == 1$}
                    \State $\beta = \texttt{true\_nrm}$ \Comment{norm of initial vector}
                    \EndIf
                    \State $V_{j:j+1} =V_{j:j+1}/ \texttt{true\_nrm}$ \Comment{scale for Arnoldi}
                    \State $\texttt{temp}_{:,2} = \texttt{temp}_{:,2}/\texttt{true\_nrm}$ \Comment{scale for Arnoldi}
                    \State $\texttt{temp}_{j:j+1,2} = \texttt{temp}_{j:j+1,2}/\texttt{true\_nrm}$
                    \If{$j > 1$} \Comment{form projection matrices}
                    \State $\texttt{temp}_{1:j-1,1} /= \texttt{true\_nrm}$ \Comment{scale for Arnoldi}
                    \State $H_{i:j,j} = \texttt{form\_T}$ \Comment{see algo. \ref{alg:Tcwy}, \ref{alg:Tncwy}, \ref{alg:Tgsmsg}}
                    \State $H_{j-1,j} = \texttt{true\_nrm}$ \Comment{set previous norm}
                    \EndIf
				\State $v_{j+1} = v_{j+1} - V_jH_{1:j,j}$ \Comment{orthogonalize}
                    \State $ss = \sum \texttt{temp}_{1:j,2}^2$ 
				\If{$\texttt{temp}_{j+1,2} < ss$}
				\State \texttt{nrm} = \texttt{MPI.AllReduce}
				\Else
				\State \texttt{nrm\_est} = $\sqrt{\texttt{temp}_{j+1,2} - ss}$ \Comment{norm estimate from \cite{hidelatency}}
				\EndIf
				\If{\texttt{nrm\_est} $\approx$ 0} \Comment{check for lucky breakdown \cite{luckybreakdowns}}
				\State happy\_breakdown = True
				\State break
				\EndIf
                    \State $v_{j+1} = v_{j+1}/\texttt{nrm\_est}$ \Comment{``normalize"}
				\EndWhile
                    \State $\texttt{nrm} = \sqrt{v_{m+1}^Tv_{m+1}}$ \Comment{\texttt{MPI\_AllReduce} norm of last vector}
                    \State $V_{m+1} = V_{m+1}/ \texttt{nrm}$
                    \State $H_{m+1,m} = \texttt{nrm}$
                
				\EndProcedure
            \end{algorithmic}\label{alg:lowsync_skeleton}
		\end{algorithm}

            \begin{algorithm}[H]
			\caption{Hybrid CWY Correction Matrix (sect. \ref{sec:cwy})}
			\begin{algorithmic}[1]
                    \Procedure{HCWY\_T}{$\texttt{temp},T$}       \Comment{inner products (line 13 algo. (\ref{alg:lowsync_skeleton})), correction matrix}
                    \State $T_{j,1:j-1} = -\texttt{temp}_{1:j-1,1}^T T_{1:j-1,1:j-1}$ \Comment{Correction matrix via (\ref{eqn:recursionForT})}
                    \State $H_{1:j,j} = T_j \texttt{temp}_{1:j,2}$ \Comment{Hessenberg values}
				\EndProcedure
			\end{algorithmic}\label{alg:Tcwy}
		\end{algorithm}
            \begin{algorithm}[H]

        \caption{Hybrid NCWY Correction Matrix (sect. \ref{sec:neumannseries}) }
        \begin{algorithmic}[1]
            \Procedure{HNCWY\_T}{$\texttt{temp}$}       \Comment{inner products (line 13 algo. (\ref{alg:lowsync_skeleton}))}
				\State $T_{j,1:j-1} = -\texttt{temp}_{1:j-1,1}^T$ \Comment{Correction matrix (\ref{subeqn:neumann_approx_t})}
                    \State $H_{1:j,j} = T_j \texttt{temp}_{1:j,2}$ \Comment{Hessenberg values}
				\EndProcedure
        \end{algorithmic}\label{alg:Tncwy}
    \end{algorithm}           

            \begin{algorithm}[H]
			\caption{Hybrid Iterated Gau{\ss}-Seidel Correction Matrix (sect. \ref{sec:gsmgs}) }
			\begin{algorithmic}[1]
				\Procedure{HGSMGS\_T}{$\texttt{temp}$}       \Comment{inner products (line 13 algo. (\ref{alg:lowsync_skeleton}))}
				\State $M_{j,1:j-1} = \texttt{temp}_{1:j-1,1}^T$
				\State $N_{1:j-1,j} = -\texttt{temp}_{1:j-1,1}$
				\State $M^{-1}_{j,1:j-1} = - \texttt{temp}_{1:j-1}^TM^{-1}_{1:j-1,1:j-1}$
                    \State $p = [I_j + N_{1:j,1:j}M^{-1}_{1:j,1:j}]\texttt{temp}_{1:j,2}$ \Comment{Part 1: mat-vec}
                    \State $H_{1:j,j} = \texttt{triang\_sol}(M_{1:j,1:j}, p)$ \Comment{Part 2: lower triangular solve + Hessenberg values}
				\EndProcedure
            \end{algorithmic}\label{alg:Tgsmsg}
		\end{algorithm}


    \section{Numerical Experiments} \label{sec:numerical_exp}





      \subsection{Application to exponential time integration} \label{sec:expint}     

        

     Over the past several decades significant amount of research on exponential integrators resulted in development of many classes of efficient and accurate methods of this type.  Initially exponential methods were used to solve diffusion or diffusion-reaction systems \cite{friesner89,gall_saad92} but with time the range of applications of these techniques increased.  Exponential integrators were used in computational weather prediction \cite{calindri2020exp,val_cube, peixoto2019semi}, combustion \cite{stewart2023variable}, astrophysics \cite{pranabmagneto, deka2023lexint, einkemmermagneto}, and other applications where the methods were shown to be competitive with state-of-the-art techniques for solving stiff PDEs. Exponential integrators of various types have been developed such as multistep \cite{val_cube}, Runge-Kutta \cite{johnthesis, valthesis, loffeld2013comparative, tokjcp}, Rosenbrock \cite{dallerit2023second, hochbruck2010exponential,  pranabmagneto, deka2023lexint}, polynomial \cite{buvoli2021exponential}, exponential parareal \cite{buvoli2024exponential, timeintegratorreview_paraexpex, huang2024parareal} and more. 

        The motivation for extending the low-synchronization methods to exponential integrators is to improve the integrator's parallel performance and strong scaling properties for solving large-scale stiff PDEs. Our primary application of interest is the shallow-water (SW) equations, but three additional examples are presented as a proof-of-concept to illustrate the benefits of the low-synchronization methods. Given that our primary application relates to weather prediction, we use exponential integrators that were developed and tested for the SW equations. The exponential multistep methods were presented in \cite{val_cube} and were shown to produce accurate, physically plausible results even with large timesteps \cite{val_cube}. In addition, we also study performance of the newly developed stiffly-resilient exponential Runge-Kutta type integrators (\texttt{SRERK}) which are constructed to mitigate issues with order reduction when solving stiff differential equations \cite{valthesis}. The new \texttt{SRERK} methods have not been studied as extensively as other exponential methods and the stiff SW equations are a good test problem for these techniques. The set of the time integrators of different order we chose is diverse enough for a thorough study of the low-synchronization methods and their effect on the parallel efficiency of exponential methods for our application. The multistep and Ruge-Kutta type integrators that are used for our experiments are listed below.   
        \begin{itemize}           
        \item $4^{th}$, $5^{th}$, and $6^{th}$ order exponential multi-step methods (\texttt{epi}) \cite{val_cube}
        \begin{subequations}\label{eqn:epi4}
               \begin{alignat}{2}
                u_{n+1} &= u_n + \varphi_1(hJ_n)hf_n + \sum_{m=1}^M \varphi_m(hJ_n)v_m \label{subeqn:advance_multi} \\
                v_m &= \sum_{i=1}^P\alpha_{m,i}hR(u_{n-i}) \\
                R(u_{n-i}) &= f_{n-i} - f_n - J_n(u_{n-i} - u_n)
               \label{subeqn:b}
            \end{alignat}
          \end{subequations} 

          with coefficients $\alpha_{m,i} $ 
            \begin{align*}
                A_4 &= \begin{pmatrix}
                    0 & 0 \\
                    \frac{-3}{10} & 
                    \frac{3}{40} \\
                    \frac{32}{5} & \frac{-11}{10}
                \end{pmatrix} \\
            A_5 &= \begin{pmatrix}
                0 & 0 & 0 \\
                -\frac{4}{5} & \frac{2}{5} & -\frac{4}{45}\\
                12 & -\frac{9}{2} & \frac{8}{9} \\
                3 & 0 & -\frac{1}{3}
            \end{pmatrix} \\
            A_6 &= \begin{pmatrix}
                0 & 0 & 0 & 0 \\
                -\frac{49}{60} & \frac{351}{560} & -\frac{359}{12600} & \frac{367}{6720} \\
                \frac{97}{7} & -\frac{99}{14} & \frac{176}{63} & -\frac{1}{2} \\
                \frac{485}{21} & -\frac{151}{14} & \frac{23}{9} & -\frac{31}{168}
            \end{pmatrix}
        \end{align*}

            \item $3^{rd}$ order Stiffly Resilient Exponential Runge-Kutta (\texttt{SRERK3}) \cite{valthesis}                       \begin{subequations}\label{eqn:epirk4}
            \begin{alignat}{2}
                z_1 &= u_n + \varphi_1\left( \frac{3}{4} h J_n \right) \frac{3}{4} h f(u_n) \\
                u_{n+1} &= u_n + \varphi_1(h J_n) h f(u_n) + \varphi_3(h J_n) \frac{32}{9}h R(z_1)
            \end{alignat}
          \end{subequations}

            \item $6^{th}$ order \texttt{SRERK6} \cite{valthesis}

            \begin{subequations}\label{eqn:epirk4}
            \begin{alignat}{2}
                z_i &= u_n + \varphi_1\left( c_i h J_n \right) c_i h f(u_n) \qquad i \in \{1,2\} \\
                z_j &= u_n + \varphi_1(c_j h J_n) c_j h f(u_n) + \varphi_3(c_j h J_n) c_j^3 \sum_{k=1}^2 \beta_{3,k}h R(z_k) \label{subeqn:srerk6_stage2_p1} \\
                &+ 
                \varphi_4(c_j h J_n) c_j^4 \sum_{k=1}^2 \beta_{4,k} h R(w_k) \qquad j \in \{3,4,5,6\} \label{subeqn:srerk6_stage2_p2}\\
                u_{n+1} &= u_n \varphi_1(h J_n) h f(u_n) + \varphi_3(h J_n) \sum_{k=3}^6 \alpha_{3,k} h R(z_k) + \varphi_4(h J_n) \sum_{k=3}^6 \alpha_{4,k} h R(z_k) \label{subeqn:srerk6_pt1}\\
                &+ \varphi_5(h J_n) \sum_{k=3}^6 \alpha_{5,k} h R(z_k) + \varphi_6(h J_n) \sum_{k=3}^6 \alpha_{6,k} h R(z_k) \label{subeqn:srerk6_pt2}
            \end{alignat}
          \end{subequations}

          with $c = [ \frac{10 - \sqrt{10}}{15}, \frac{10 + \sqrt{10}}{15}, \frac{1}{4}, \frac{1}{2}, \frac{2}{4}, 1, ]$, 

          \begin{align*}
              \beta_{i,k} = \begin{bmatrix}
                  \frac{155 + 65\sqrt{10}}{18} & \frac{155 - 65\sqrt{10}}{18} \\
                  \frac{-100 - 55\sqrt{10}}{4} & \frac{-100 - 55\sqrt{10}}{4}
              \end{bmatrix} \qquad \alpha_{j,k} = \begin{bmatrix}
                  128 & -48 & \frac{128}{9} & -2 \\
                  -1664 & 912 & -\frac{896}{3} & 44 \\
                  9216 & -6144 & \frac{7168}{3} & -384 \\
                  -20480 & 15360 & -\frac{20480}{3} & 1280
              \end{bmatrix}
          \end{align*}
          
            \end{itemize}

        Both the multistep and Runge-Kutta -type integrators require the computation of a linear combination of the exponential-like $\varphi$-functions times a vector (sec. \ref{sec:linearPhiFunc}). To compute linear combinations (\ref{eq:ansatz}), our algorithm of choice is \texttt{KIOPS} (outlined in Algorithm \ref{alg:kiops}, referenced in \ref{sec:io}) \cite{kiops}. The integrators are implemented in order to minimize the number of \texttt{KIOPS} calls per timestep. Each multistep method will require one call to \texttt{KIOPS} to compute the linear combination of $\varphi$-functions in order to advance the solution to the next timestep $t_{n+1}$ (e.g. equation (\ref{subeqn:advance_multi}) ). \texttt{SRERK3} will require two calls to \texttt{KIOPS}, one for the stage and one to advance the solution. Finally, \texttt{SRERK6} will require six calls to \texttt{KIOPS}, one to advance the solution (eg. (\ref{subeqn:srerk6_pt1}) - (\ref{subeqn:srerk6_pt2}) ), one for the stages $z_i$, and four calls to stage $z_j$ (e.g. (\ref{subeqn:srerk6_stage2_p1}) - (\ref{subeqn:srerk6_stage2_p2}) ). 
                
        To compare the performance of the low-synchronization GS variants, we implement the orthogonalization methods discussed in sections \ref{sec:cwy} to \ref{sec:gsmgs} within the \texttt{KIOPS} algorithm.  A parallel implementation based on an incomplete-orthogonalization Classical Gram-Schmidt (\texttt{ioCGS}) is provided in the supplementary material of \cite{val_cube}. \texttt{KIOPS} is adapted in line 3 of Algorithm \ref{alg:kiops} to compare the efficiency of exponential integrators with the full-orthogonalization, low-synchronization, MGS variants against the reference \texttt{ioCGS} algorithm. The adaptivity for $\tau_i$ and $m_i$ remains the same throughout the experiments, allowing for a fair comparison between the orthogonalization techniques.

        \begin{algorithm}[H]
        \caption{Krylov Incomplete Orthogonalization Process (\texttt{KIOPS}) \cite{kiops}}
        \begin{algorithmic}[1]
            \Procedure{KIOPS}{$A,b, m, \tau_{final}$}       \Comment{matrix, initial vector, size of subspace, final $\tau$}
                \While{substeping in $\tau$}
            \State Compute Arnoldi Decomposition:  $\texttt{Arnoldi}(A, b,m)$ (eg. alg. (\ref{al:mgs}) )
            \State Compute matrix exponential using (\ref{eq:kryexp})
            \State Compute local error estimate for step $\tau_i$ (eqn. (36) in \cite{kiops})
            \State Compute suggested new parameters $\tau_{new}$ and $m_{new}$ (eqn. (44) and (46) in \cite{kiops})
                \State Decide to accept or reject solution at step $\tau_i$ (sect. 3.5 in \cite{kiops})
            \EndWhile			
            \EndProcedure
        \end{algorithmic}\label{alg:kiops}
    \end{algorithm}

    \subsection{Experimental set-up}

    We highlight the strong-scaling parallel efficiency of the five exponential integrators with each of the hybrid low-synchronization methods presented in sections \ref{sec:cwy} to \ref{sec:gsmgs}. To complement the strong scaling\footnote{Strong scaling refers to the application's performance or runtime when the problem size is fixed and the number of processors varies. } plots, the runtime ratio between the reference algorithm \texttt{ioCGS} and each of the hybrid low-synchronization methods is presented in a table. The runtime ratio for a low-synchronization method \textit{ls} is computed as 
    \[\texttt{runtime\_ratio}_{ls} = \texttt{runtime}_{ls}/ \texttt{runtime}_{ioCGS} \] 
    where a runtime ratio below $1.0$ implies that the exponential integrator with the low-synchronization method is more efficient than a version with an incomplete Classical Gram-Schmidt orthogonalization scheme.  
    To avoid pollution of the data from the machine state variability, 5 runs are performed and the median is reported in strong scaling plots for each problem and integrator.  

    Performance studies were conducted on the Lenovo ThinkSystem supercomputer of Shared Services Canada. This cluster runs Linux on Intel® Xeon® Platinum 8380 CPU 2.30GHz processors, equipped with 512 GB of memory per node, and utilizes an Infiniband HDR network.

    \subsection{Moderate Scale Example PDEs}\label{sec:expdes}

    For the proof-of-concept tests of the performance of the low-synchronization methods compared to standard Krylov algorithms, we choose problems that are routinely used to study the performance of stiff integrators: 
    \begin{itemize}
        \item Allen-Cahn (AC) equation \cite{acref}
            \begin{equation}
        \frac{\partial u}{\partial t} = \epsilon \nabla^2 u + u - u^3
    \end{equation} with no-flow boundary conditions on $x,y \in [-1,1]$ and a diffusion coefficient of $\epsilon = 0.1$. The initial condition is given by 
    \begin{equation}
        u(x,y,t-0) = 0.1 + 0.1\cos(2 \pi x) \cos(2 \pi y)
    \end{equation}
    with a final time of integration to be $t_f = 0.02$ and a timestep of $dt = 10^{-2}$. 
    
        \item Advection-Diffusion-Reaction (ADR) \cite{adrref}
            \begin{equation}
        \frac{\partial u}{\partial t} = \epsilon \nabla^2 u - \alpha(u_x + u_y) + \gamma u(u-\frac{1}{2})(1-u)
    \end{equation}
    where homogeneous Neumann boundary conditions are imposed on $x, y \in [0,1]$. The coefficients $\epsilon$ and $\alpha$ measure the diffusion and advection respectively. For this scenario, we consider an advection-dominated regime with $\epsilon = 10^{-2}$ and $\alpha = -6 $, with the initial condition given by
    \begin{equation}
        u(x,y,t=0) = 0.3 + 256(xy(1-x)(1-y))^2.
    \end{equation}
    We chose the final time of integration to be $t_f = 10^{-1}$ with a timestep of $dt = 10^{-2}$
        
        \item Burger's equation (BURG)
            \begin{equation}
        \frac{\partial u}{\partial t} = - \frac{1}{2} \left( \frac{\partial u^2}{\partial x} + \frac{\partial u^2}{\partial y}  \right) + \epsilon \nabla^2 u
    \end{equation} with Dirichlet boundary conditions on $x,y \in [0,1]$. The initial condition is given by 
    \begin{equation}
        u(x,y,t=0) = \sin(3 \pi x)^2\sin(3 \pi y)^2 (1-y).
    \end{equation} We choose a final time of $t_f = 0.04$ with a timestep of $dt = 10^{-2}$, and $\epsilon = 10^{-3}$.
    
    \end{itemize}
    In all problems, the number of degrees of freedom is $N = 2000^2$, and second order finite differences are used to approximate the first and second spatial derivatives. Strong scaling results with processor sizes $p = [100, 400, 625, 1600, 2500]$ are presented, and a Krylov tolerance of $10^{-12}$ is used within \texttt{KIOPS}.

    \subsubsection{Strong Scaling Results}

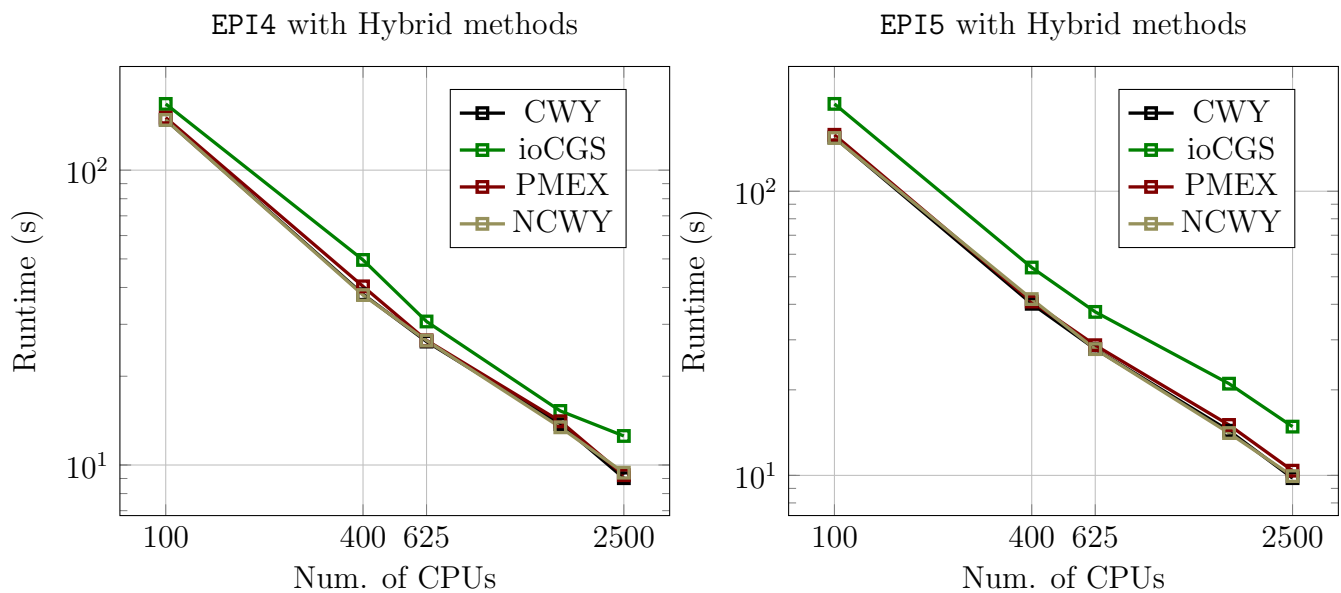
\begin{figure}[H]
          \centering
          \begin{subfigure}[b]{0.5\textwidth}
          \begin{tikzpicture}
        \begin{axis}[
        grid=major, width=1\textwidth, height=0.85\textwidth,
        title={\texttt{EPI4} with Hybrid methods},
        xtick={100, 400, 625, 2500},
        xticklabels={100, 400, 625, 2500},
        xlabel={Num. of CPUs}, ylabel={Runtime (s)},
        xmode=log, ymode=log,
        legend style={at={(0.60,0.95)},anchor=north west}
        ]
        \addplot[color=black!80!black, line width = 1.25pt, mark=square]
        table{ac/epi4/med_cwyne1s.txt};
        \addlegendentry{CWY};
        \addplot[color=green!50!black, line width = 1.25pt, mark=square]
        table{ac/epi4/med_kiops.txt};
        \addlegendentry{ioCGS};
        \addplot[color=red!50!black, line width = 1.25pt, mark=square]table{ac/epi4/med_pmexne1s.txt};
        \addlegendentry{PMEX};
        \addplot[color=yellow!50!black, line width = 1.25pt, mark=square]table{ac/epi4/med_icwyne1s.txt};
        \addlegendentry{NCWY};
        \end{axis}
        \end{tikzpicture}
          \end{subfigure}%
         \hfill
          \begin{subfigure}[b]{0.5\textwidth}
              \begin{tikzpicture}
        \begin{axis}[
        grid=major, width=1\textwidth, height=0.85\textwidth,
        xtick={100, 400, 625, 2500},
        xticklabels={100, 400, 625, 2500},
        xlabel={Num. of CPUs}, ylabel={Runtime (s)},
        title={\texttt{EPI5} with Hybrid methods},
        xmode=log, ymode=log,
        legend style={at={(0.60,0.95)},anchor=north west}
        ]
        \addplot[color=black!80!black, line width = 1.25pt, mark=square]
        table{ac/epi5/med_cwyne1s.txt};
        \addlegendentry{CWY};
        \addplot[color=green!50!black, line width = 1.25pt, mark=square]
        table{ac/epi5/med_kiops.txt};
        \addlegendentry{ioCGS};
        \addplot[color=red!50!black, line width = 1.25pt, mark=square]table{ac/epi5/med_pmexne1s.txt};
        \addlegendentry{PMEX};
        \addplot[color=yellow!50!black, line width = 1.25pt, mark=square]table{ac/epi5/med_icwyne1s.txt};
        \addlegendentry{NCWY};
        \end{axis}
        \end{tikzpicture}
          \end{subfigure}
          \caption{Strong scaling results for AC with low-synchronization hybrid methods for \texttt{epi4} and \texttt{epi5}. }
          \label{fig:ac_epi45}
        \end{figure}

\begin{figure}[H]
          \centering
          \begin{subfigure}[b]{0.5\textwidth}
           \begin{tikzpicture}
        \begin{axis}[
        grid=major, width=1\textwidth, height=0.85\textwidth,
        xtick={100, 400, 625, 2500},
        xticklabels={100, 400, 625, 2500},
        xlabel={Num. of CPUs}, ylabel={Runtime (s)},
        title={\texttt{srerk3} with Hybrid methods},
        xmode=log, ymode=log,
        legend style={at={(0.60,0.95)},anchor=north west}
        ]
        \addplot[color=black!80!black, line width = 1.25pt, mark=square]
        table{ac/srerk3/med_cwyne1s.txt};
        \addlegendentry{CWY};
        \addplot[color=green!50!black, line width = 1.25pt, mark=square]
        table{ac/srerk3/med_kiops.txt};
        \addlegendentry{ioCGS};
        \addplot[color=red!50!black, line width = 1.25pt, mark=square]table{ac/srerk3/med_pmexne1s.txt};
        \addlegendentry{GSMGS};
        \addplot[color=yellow!50!black, line width = 1.25pt, mark=square]table{ac/srerk3/med_icwyne1s.txt};
        \addlegendentry{NCWY};
        \end{axis}
        \end{tikzpicture}
          \end{subfigure}%
         \hfill
          \begin{subfigure}[b]{0.5\textwidth}
              \begin{tikzpicture}
        \begin{axis}[
        grid=major, width=1\textwidth, height=0.85\textwidth,
        xtick={100, 400, 625, 2500},
        xticklabels={100, 400, 625, 2500},
        xlabel={Num. of CPUs}, ylabel={Runtime (s)},
        title={\texttt{srerk6} with Hybrid methods},
        xmode=log, ymode=log,
        legend style={at={(0.60,0.95)},anchor=north west}
        ]
        \addplot[color=black!80!black, line width = 1.25pt, mark=square]
        table{ac/srerk6/med_cwyne1s.txt};
        \addlegendentry{CWY};
        \addplot[color=green!50!black, line width = 1.25pt, mark=square]
        table{ac/srerk6/med_kiops.txt};
        \addlegendentry{ioCGS};
        \addplot[color=red!50!black, line width = 1.25pt, mark=square]table{ac/srerk6/med_pmexne1s.txt};
        \addlegendentry{GSMGS};
        \addplot[color=yellow!50!black, line width = 1.25pt, mark=square]table{ac/srerk6/med_icwyne1s.txt};
        \addlegendentry{NCWY};
        \end{axis}
        \end{tikzpicture}
          \end{subfigure}
          \caption{Strong scaling results for AC with low-synchronization hybrid methods for \texttt{srerk} integrators.}
          \label{fig:ac_srerk}
        \end{figure}
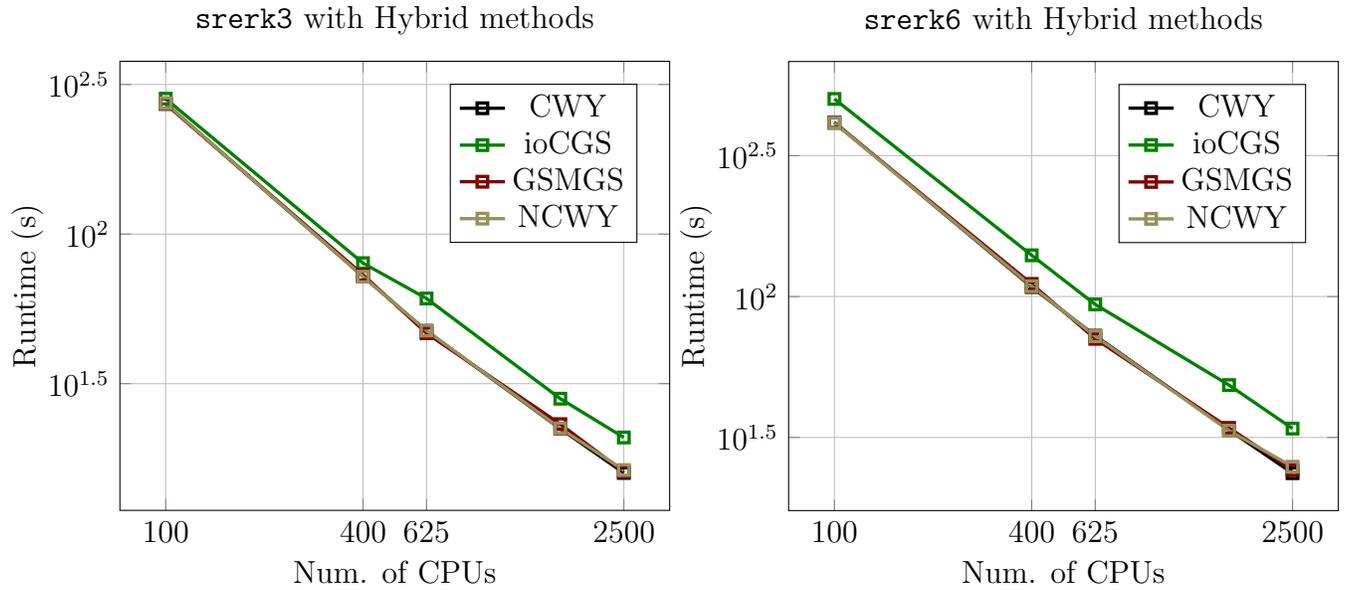

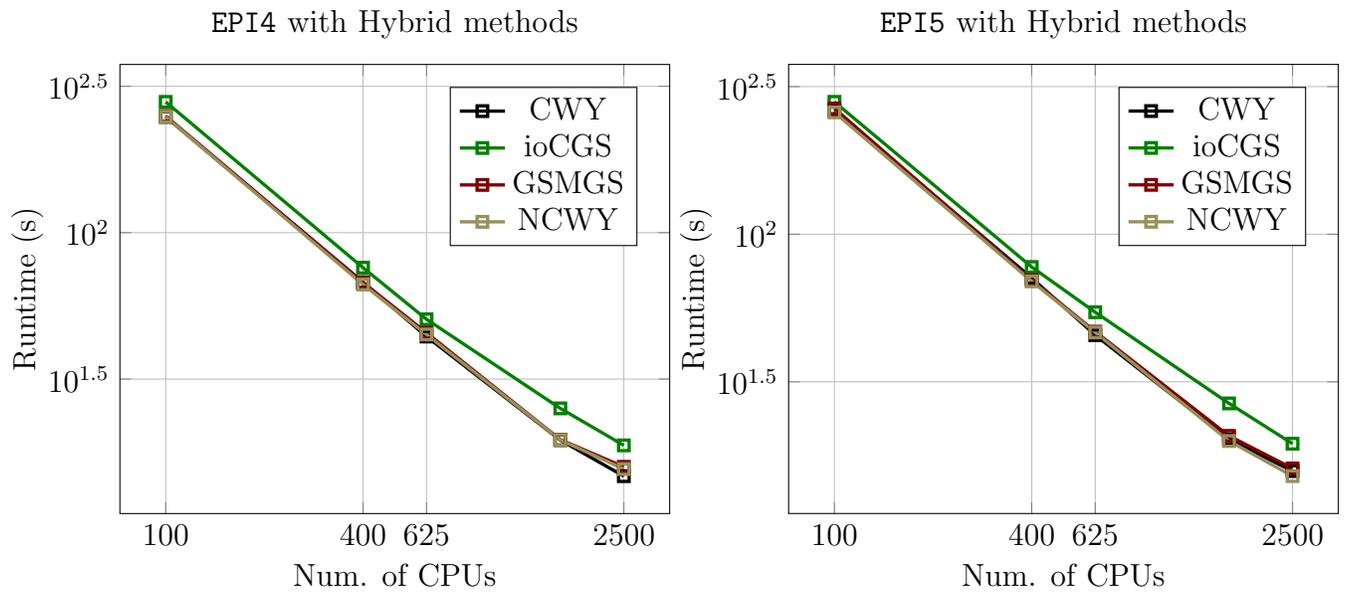
\begin{figure}[H]
          \centering
          \begin{subfigure}[b]{0.5\textwidth}
          \begin{tikzpicture}
        \begin{axis}[
        grid=major, width=1\textwidth, height=0.85\textwidth,
        title={\texttt{EPI4} with Hybrid methods},
        xtick={100, 400, 625,  2500},
        xticklabels={100, 400, 625, 2500},
        xlabel={Num. of CPUs}, ylabel={Runtime (s)},
        xmode=log, ymode=log,
        legend style={at={(0.60,0.95)},anchor=north west}
        ]
        \addplot[color=black!80!black, line width = 1.25pt, mark=square]
        table{adr/epi4/med_cwyne1s.txt};
        \addlegendentry{CWY};
        \addplot[color=green!50!black, line width = 1.25pt, mark=square]
        table{adr/epi4/med_kiops.txt};
        \addlegendentry{ioCGS};
        \addplot[color=red!50!black, line width = 1.25pt, mark=square]table{adr/epi4/med_pmexne1s.txt};
        \addlegendentry{GSMGS};
        \addplot[color=yellow!50!black, line width = 1.25pt, mark=square]table{adr/epi4/med_icwyne1s.txt};
        \addlegendentry{NCWY};
        \end{axis}
        \end{tikzpicture}
          \end{subfigure}%
         \hfill
          \begin{subfigure}[b]{0.5\textwidth}
              \begin{tikzpicture}
        \begin{axis}[
        grid=major, width=1\textwidth, height=0.85\textwidth,
        xtick={100, 400, 625, 2500},
        xticklabels={100, 400, 625, 2500},
        xlabel={Num. of CPUs}, ylabel={Runtime (s)},
        title={\texttt{EPI5} with Hybrid methods},
        xmode=log, ymode=log,
        legend style={at={(0.60,0.95)},anchor=north west}
        ]
        \addplot[color=black!80!black, line width = 1.25pt, mark=square]
        table{adr/epi5/med_cwyne1s.txt};
        \addlegendentry{CWY};
        \addplot[color=green!50!black, line width = 1.25pt, mark=square]
        table{adr/epi5/med_kiops.txt};
        \addlegendentry{ioCGS};
        \addplot[color=red!50!black, line width = 1.25pt, mark=square]table{adr/epi5/med_pmexne1s.txt};
        \addlegendentry{GSMGS};
        \addplot[color=yellow!50!black, line width = 1.25pt, mark=square]table{adr/epi5/med_icwyne1s.txt};
        \addlegendentry{NCWY};
        \end{axis}
        \end{tikzpicture}
          \end{subfigure}
          \caption{Strong scaling results for ADR with low-synchronization hybrid methods for \texttt{epi4} and \texttt{epi5}.}
          \label{fig:adr_epi45}
        \end{figure}

\begin{figure}[H]
          \centering
          \begin{subfigure}[b]{0.5\textwidth}
           \begin{tikzpicture}
        \begin{axis}[
        grid=major, width=1\textwidth, height=0.85\textwidth,
        xtick={100, 400, 625, 2500},
        xticklabels={100, 400, 625, 2500},
        xlabel={Num. of CPUs}, ylabel={Runtime (s)},
        title={\texttt{srerk3} with Hybrid methods},
        xmode=log, ymode=log,
        legend style={at={(0.60,0.95)},anchor=north west}
        ]
        \addplot[color=black!80!black, line width = 1.25pt, mark=square]
        table{adr/srerk3/med_cwyne1s.txt};
        \addlegendentry{CWY};
        \addplot[color=green!50!black, line width = 1.25pt, mark=square]
        table{adr/srerk3/med_kiops.txt};
        \addlegendentry{ioCGS};
        \addplot[color=red!50!black, line width = 1.25pt, mark=square]table{adr/srerk3/med_pmexne1s.txt};
        \addlegendentry{GSMGS};
        \addplot[color=yellow!50!black, line width = 1.25pt, mark=square]table{adr/srerk3/med_icwyne1s.txt};
        \addlegendentry{NCWY};
        \end{axis}
        \end{tikzpicture}
          \end{subfigure}%
         \hfill
          \begin{subfigure}[b]{0.5\textwidth}
              \begin{tikzpicture}
        \begin{axis}[
        grid=major, width=1\textwidth, height=0.85\textwidth,
        xtick={100, 400, 625, 2500},
        xticklabels={100, 400, 625, 2500},
        xlabel={Num. of CPUs}, ylabel={Runtime (s)},
        title={\texttt{srerk6} with Hybrid methods},
        xmode=log, ymode=log,
        legend style={at={(0.60,0.95)},anchor=north west}
        ]
        \addplot[color=black!80!black, line width = 1.25pt, mark=square]
        table{adr/srerk6/med_cwyne1s.txt};
        \addlegendentry{CWY};
        \addplot[color=green!50!black, line width = 1.25pt, mark=square]
        table{adr/srerk6/med_kiops.txt};
        \addlegendentry{ioCGS};
        \addplot[color=red!50!black, line width = 1.25pt, mark=square]table{adr/srerk6/med_pmexne1s.txt};
        \addlegendentry{GSMGS};
        \addplot[color=yellow!50!black, line width = 1.25pt, mark=square]table{adr/srerk6/med_icwyne1s.txt};
        \addlegendentry{NCWY};
        \end{axis}
        \end{tikzpicture}
          \end{subfigure}
          \caption{Strong scaling results for ADR with low-synchronization hybrid methods for \texttt{srerk3} and \texttt{srerk6}.}
          \label{fig:adr_srerk}
        \end{figure}
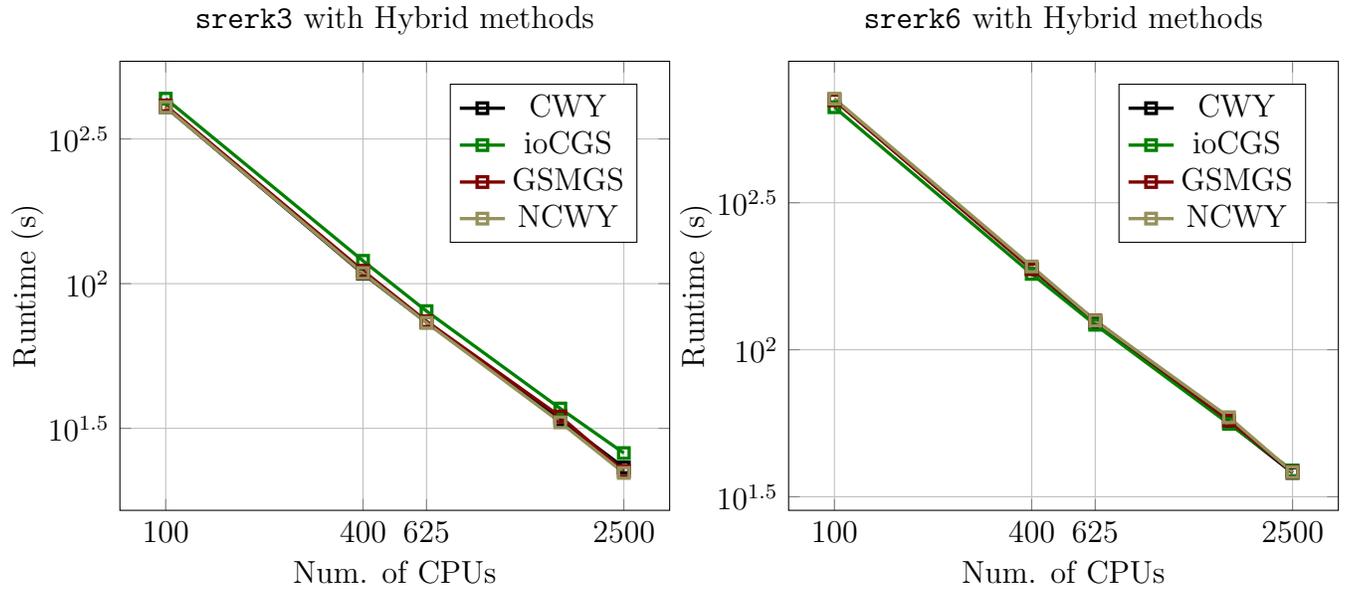

    The strong scaling results for AC are shown in Figures \ref{fig:ac_epi45} - \ref{fig:ac_srerk}, for ADR in \ref{fig:adr_epi45} - \ref{fig:adr_srerk}, and BURG in Figures \ref{fig:burg_epi45} - \ref{fig:burg_srerk}, with the tables for runtime rations displayed in Tables \ref{runtimetratioac}, \ref{runtimetratioadr}, and \ref{runtimetratioburg} respectively. 

    Overall, the degree of improvement varied for each example problem, however, it is clear that the low-synchronization methods perform on-par if not better than the reference incomplete CGS algorithm. Our numerical experiments show that the greater computational complexity of the low-synchronization algorithms does not significantly slow down the runtime for these moderately sized example problems. The AC example displays the largest amount of runtime savings for all integrators. Meanwhile, the runtime for Burger's equation with the low-synchronization algorithms remains nearly identical to the reference \texttt{ioCGS}. This is in-line with our performance expectations for these moderately sized problems on a relatively small number of processors, given that the low-synchronization methods were designed for massively parallel systems and large scale problems. It is informative to see that there is no penalty for using these methods on a purely distributed memory framework; we are not using GPUs. Now we discuss the results of each individual problem in greater detail. 
    

    For AC, both the \texttt{epi} and \texttt{srerk} integrators showed significant improvement with all three low-synchronization methods (see Figures \ref{fig:ac_epi45} - \ref{fig:ac_srerk}). The greatest improvement for an \texttt{epi} integrator is seen for \texttt{epi6} with a $40\%$ reduction in runtime with \texttt{H-NCWY} with $2500$ processors (see Table \ref{runtimetratioac}). For the \texttt{srerk} integrators (Figure \ref{fig:ac_srerk}), the greatest improvement is seen for \texttt{srerk6} with a $31\%$ reduction in runtime with \texttt{H-CWY} with $1600$ and $2500$ processors. This test problem is a clear example of the improved performance of exponential integrators with a low-synchronization Arnoldi method.
    
    For ADR, there are improvements with the \texttt{epi} integrators (see Figure \ref{fig:adr_epi45} - \ref{fig:adr_e6}) and marginal to no improvement with the \texttt{srerk} integrators (see Figure \ref{fig:adr_srerk}). The greatest improvement for an \texttt{epi} integrator is seen for \texttt{epi5} with a $25\%$ reduction in runtime with \texttt{H-NCWY} with $1600$ processors (see Table \ref{runtimetratioadr}). For the \texttt{srerk} integrators, there were marginal improvements with \texttt{srerk3}, at most $15\%$ reduction in runtime with \texttt{H-NCWY} for $2500$ processors (see Figure \ref{fig:adr_srerk} and Table \ref{runtimetratioadr}). Meanwhile, the results for \texttt{srerk6} displayed little to no improvement. For processor sizes less than $2500$, all methods were slightly slower compared to \texttt{ioCGS} and showed marginal improvement for for $2500$ processors.  
    
    For Burger's equation, both \texttt{epi} and \texttt{srerk} integrators performed similarly between the low-synchronization method and \texttt{ioCGS}. The most significant improvement over all integrators is shown for \texttt{epi4} with $2500$ processors with a $25\%$ reduction in runtime with \texttt{H-NCWY}. For the \texttt{srerk} integrators, the strong scaling curves display nearly overlapping lines, indicating the runtime of the low-synchronzation methods are on-par with the reference \texttt{ioCGS}. Small improvements ($\leq 9\%$) are seen for both \texttt{srerk} integrators with \texttt{H-NCWY} and \texttt{H-CWY} for processors larger than $100$. 
    
    For this tests set, it remains unclear which low-synchronization methods is more favorable. In the AC and ADR \texttt{epi} strong scaling results (Figures \ref{fig:ac_epi45} - \ref{fig:ac_e6} and \ref{fig:adr_epi45} - \ref{fig:adr_e6}) (where the three low-synchronization methods showed to be more efficient) the runtime between the three low-synchronization methods is nearly identical, as shown by the strong scaling curves are nearly overlapping. However, for the largest processor size the best results are achieved with either \texttt{H-CWY} or \texttt{H-NCWY}, (see Tables \ref{runtimetratioac} and \ref{runtimetratioadr}).

    
        \begin{figure}[H]
          \centering
          \begin{tikzpicture}
        \begin{axis}[
        grid=major, width=0.52\textwidth, height=0.43\textwidth,
        title={\texttt{EPI6} with Hybrid methods},
        xtick={100, 400, 625, 2500},
        xticklabels={100, 400, 625, 2500},
        xlabel={Num. of CPUs}, ylabel={Runtime (s)},
        xmode=log, ymode=log,
        legend style={at={(0.60,0.95)},anchor=north west}
        ]
        \addplot[color=black!80!black, line width = 1.25pt, mark=square]
        table{ac/epi6/med_cwyne1s.txt};
        \addlegendentry{CWY};
        \addplot[color=green!50!black, line width = 1.25pt, mark=square]
        table{ac/epi6/med_kiops.txt};
        \addlegendentry{ioCGS};
        \addplot[color=red!50!black, line width = 1.25pt, mark=square]table{ac/epi6/med_pmexne1s.txt};
        \addlegendentry{GSMGS};
        \addplot[color=yellow!50!black, line width = 1.25pt, mark=square]table{ac/epi6/med_icwyne1s.txt};
        \addlegendentry{NCWY};
        \end{axis}
        \end{tikzpicture}
          \caption{Strong scaling results for AC with low-synchronization hybrid methods for \texttt{epi6}.}
          \label{fig:ac_e6}
        \end{figure}
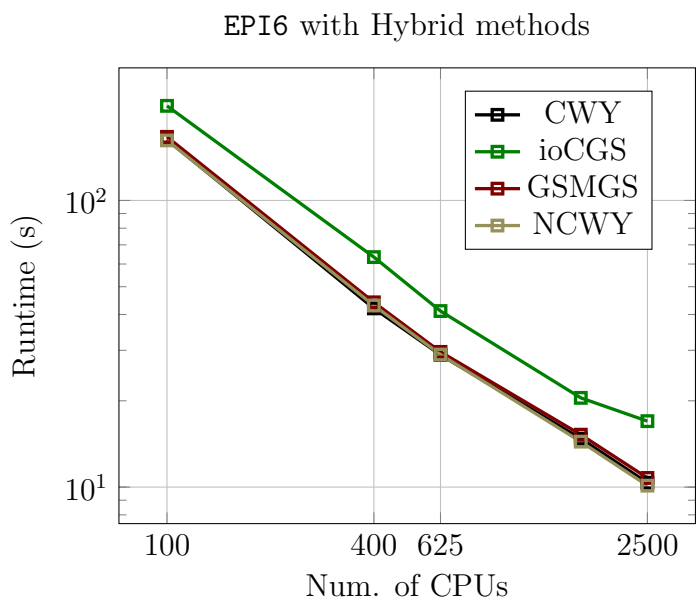
 
\begin{table}[H]
	\centering
	\begin{tabular}{c|c|c|c|c|c|c p{1cm}p{1cm}p{1cm}p{1cm}p{1cm}p{1cm}p{1cm}|}
		\toprule
		low-synchronization Method & $p$ & EPI4 & EPI5 & EPI6 & SRERK3 & SRERK6  \\ \hline
		\multirow{ 2}{*}{\texttt{H-GSMGS}} 
		& 100 & 0.90 & 0.78 & 0.78 & 0.96 & 0.82  \\
		& 400 & 0.81 & 0.76  & 0.70 & 0.91 & 0.79  \\
		& 625 & 0.86 & 0.76 & 0.72 & 0.77 & 0.75\\
		& 1600 & 0.92 & 0.72 & 0.74 & 0.82 & 0.71\\
		& 2500 & 0.74 & 0.70 & 0.63 & 0.77 & 0.71\\
		\midrule
		\multirow{ 6}{*}{\texttt{H-NCWY}} 
		& 100 & 0.88 & 0.76 & 0.76 & 0.97 & 0.82 \\
		& 400 & 0.76 & 0.77  & 0.68 & 0.90 & 0.78 \\ 
		& 625 & 0.86 & 0.74 & 0.70 &  0.78 & 0.77 \\
		& 1600 & 0.88 & 0.67 & 0.70 & 0.80 & 0.69 \\
		& 2500 & 0.75 & 0.67 & 0.60 & 0.78 & 0.73 \\
		\midrule
		\multirow{ 6}{*}{\texttt{H-CWY}} 
		& 100 & 0.88 & 0.76 &  0.76 & 0.96 & 0.83  \\
		& 400 & 0.77 & 0.75  & 0.66 & 0.90 & 0.77 \\
		& 625 & 0.85 & 0.74 & 0.70 & 0.78 & 0.77\\
		& 1600 & 0.90 & 0.68 & 0.72 & 0.79 & 0.69 \\
		& 2500 & 0.72 & 0.66 & 0.61 & 0.76 & 0.69 \\
		\midrule 
	\end{tabular}
	\caption{Runtime ratios of \texttt{low-synchronization}/\texttt{ioCGS} for Allen-Cahn pde.}
	\label{runtimetratioac}
\end{table}

        \begin{figure}[H]
          \centering
          \begin{tikzpicture}
        \begin{axis}[
        grid=major, width=0.52\textwidth, height=0.43\textwidth,
        title={\texttt{EPI6} with Hybrid methods},
        xtick={100, 400, 625, 2500},
        xticklabels={100, 400, 625, 2500},
        xlabel={Num. of CPUs}, ylabel={Runtime (s)},
        xmode=log, ymode=log,
        legend style={at={(0.60,0.95)},anchor=north west}
        ]
        \addplot[color=black!80!black, line width = 1.25pt, mark=square]
        table{adr/epi6/med_cwyne1s.txt};
        \addlegendentry{CWY};
        \addplot[color=green!50!black, line width = 1.25pt, mark=square]
        table{adr/epi6/med_kiops.txt};
        \addlegendentry{ioCGS};
        \addplot[color=red!50!black, line width = 1.25pt, mark=square]table{adr/epi6/med_pmexne1s.txt};
        \addlegendentry{GSMGS};
        \addplot[color=yellow!50!black, line width = 1.25pt, mark=square]table{adr/epi6/med_icwyne1s.txt};
        \addlegendentry{NCWY};
        \end{axis}
        \end{tikzpicture}
          \caption{Strong scaling results for ADR with low-synchronization hybrid methods for \texttt{epi6}.}
          \label{fig:adr_e6}
        \end{figure}
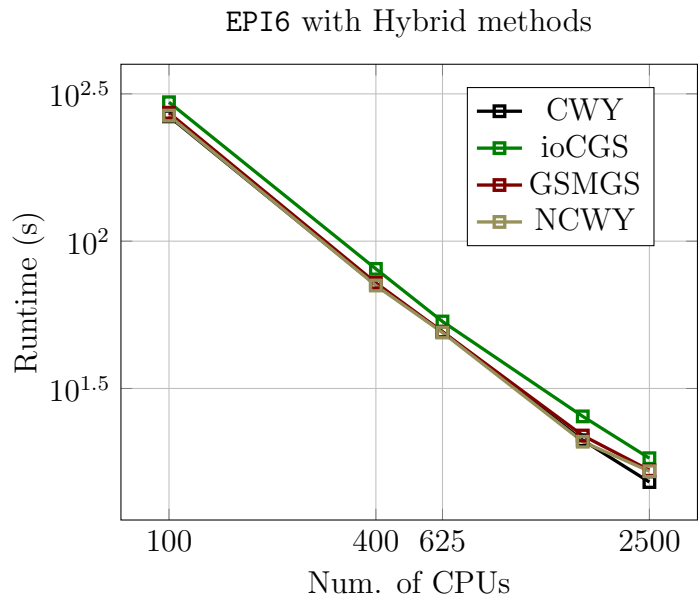

        \begin{table}[H]
	\centering
	\begin{tabular}{c|c|c|c|c|c|c p{1cm}p{1cm}p{1cm}p{1cm}p{1cm}p{1cm}p{1cm}|}
		\toprule
		low-synchronization Method & $p$ & EPI4 & EPI5 & EPI6 & SRERK3 & SRERK6  \\ \hline
		\multirow{ 2}{*}{\texttt{H-GSMGS}} 
		& 100 & 0.90 & 0.95 & 0.92 &  0.95 & 1.06   \\
		& 400 & 0.89 & 0.90  & 0.90 & 0.92 & 1.04   \\
		& 625 & 0.90 & 0.86 & 0.92 & 0.92 & 1.03  \\
		& 1600 & 0.78 & 0.78 & 0.86 & 0.94 & 1.03  \\
		& 2500 & 0.84 & 0.82 & 0.91 & 0.87 & 0.99  \\
		\hline
		\multirow{ 6}{*}{\texttt{H-NCWY}} 
		& 100 & 0.89 & 0.92 & 0.90 & 0.94 & 1.07   \\
		& 400 & 0.88 & 0.89  & 0.88 & 0.91 & 1.06    \\ 
		& 625 & 0.89 & 0.85 & 0.92 & 0.91 & 1.03  \\
		& 1600 & 0.78 & 0.75 & 0.82 & 0.90 & 1.05  \\
		& 2500 & 0.83 & 0.78 & 0.90 & 0.85 & 0.99\\
		\midrule
		\multirow{ 6}{*}{\texttt{H-CWY}} 
		& 100 & 0.89 & 0.93 &  0.90 & 0.94 & 1.06   \\
		& 400 & 0.89 & 0.92  & 0.88 & 0.90 & 1.04    \\
		& 625 & 0.87 & 0.84 & 0.92 & 0.91 & 1.08  \\
		& 1600 & 0.78 & 0.76 & 0.83 & 0.91 &1.0  \\
		& 2500 & 0.78 & 0.81 & 0.83 &0.89 & 0.99   \\
		\midrule 
	\end{tabular}
	\caption{Runtime ratios of \texttt{low-synchronization}/\texttt{ioCGS} for Advection-Diffusion-Reaction.}
	\label{runtimetratioadr}
\end{table}


        \begin{figure}[H]
          \centering
          \begin{subfigure}[b]{0.5\textwidth}
          \begin{tikzpicture}
        \begin{axis}[
        grid=major, width=1\textwidth, height=0.85\textwidth,
        title={\texttt{EPI4} with Hybrid methods},
        xtick={100, 400, 625, 2500},
        xticklabels={100, 400, 625, 2500},
        xlabel={Num. of CPUs}, ylabel={Runtime (s)},
        xmode=log, ymode=log,
        legend style={at={(0.60,0.95)},anchor=north west}
        ]
        \addplot[color=black!80!black, line width = 1.25pt, mark=square]
        table{burg/epi4/med_cwyne1s.txt};
        \addlegendentry{CWY};
        \addplot[color=green!50!black, line width = 1.25pt, mark=square]
        table{burg/epi4/med_kiops.txt};
        \addlegendentry{ioCGS};
        \addplot[color=red!50!black, line width = 1.25pt, mark=square]table{burg/epi4/med_pmexne1s.txt};
        \addlegendentry{GSMGS};
        \addplot[color=yellow!50!black, line width = 1.25pt, mark=square]table{burg/epi4/med_icwyne1s.txt};
        \addlegendentry{NCWY};
        \end{axis}
        \end{tikzpicture}
          \end{subfigure}%
         \hfill
          \begin{subfigure}[b]{0.5\textwidth}
              \begin{tikzpicture}
        \begin{axis}[
        grid=major, width=1\textwidth, height=0.85\textwidth,
        xtick={100, 400, 625, 2500},
        xticklabels={100, 400, 625, 2500},
        xlabel={Num. of CPUs}, ylabel={Runtime (s)},
        title={\texttt{EPI5} with Hybrid methods},
        xmode=log, ymode=log,
        legend style={at={(0.60,0.95)},anchor=north west}
        ]
        \addplot[color=black!80!black, line width = 1.25pt, mark=square]
        table{burg/epi5/med_cwyne1s.txt};
        \addlegendentry{CWY};
        \addplot[color=green!50!black, line width = 1.25pt, mark=square]
        table{burg/epi5/med_kiops.txt};
        \addlegendentry{ioCGS};
        \addplot[color=red!50!black, line width = 1.25pt, mark=square]table{burg/epi5/med_pmexne1s.txt};
        \addlegendentry{GSMGS};
        \addplot[color=yellow!50!black, line width = 1.25pt, mark=square]table{burg/epi5/med_icwyne1s.txt};
        \addlegendentry{NCWY};
        \end{axis}
        \end{tikzpicture}
          \end{subfigure}
          \caption{Strong scaling results for Burger's with low-synchronization hybrid methods for \texttt{epi4} and \texttt{epi5}.}
          \label{fig:burg_epi45}
        \end{figure}
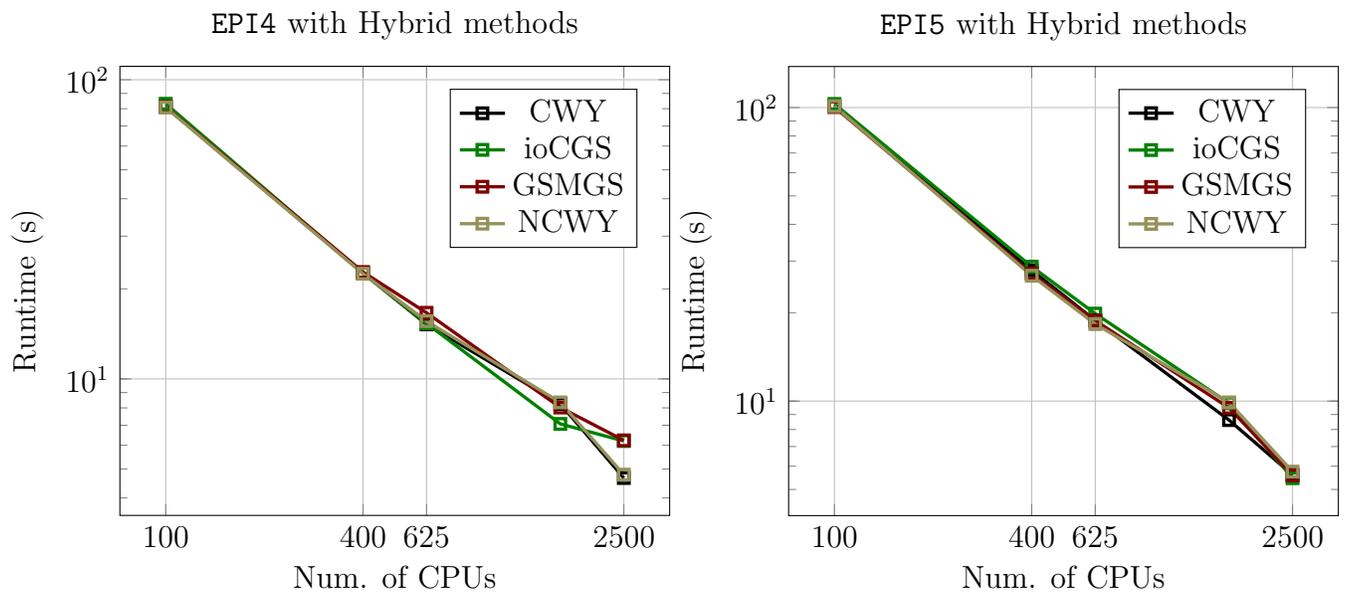
        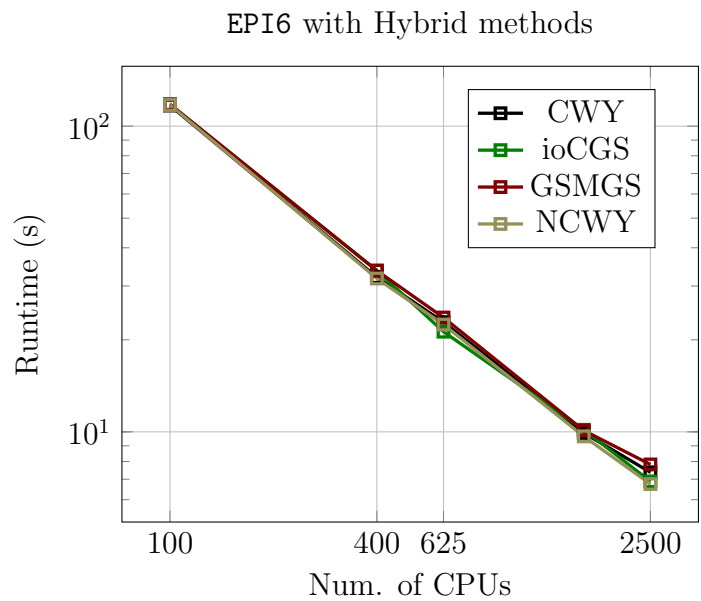
\begin{figure}[H]
          \centering
          \begin{tikzpicture}
        \begin{axis}[
        grid=major, width=0.52\textwidth, height=0.43\textwidth,
        title={\texttt{EPI6} with Hybrid methods},
        xtick={100, 400, 625, 2500},
        xticklabels={100, 400, 625, 2500},
        xlabel={Num. of CPUs}, ylabel={Runtime (s)},
        xmode=log, ymode=log,
        legend style={at={(0.60,0.95)},anchor=north west}
        ]
        \addplot[color=black!80!black, line width = 1.25pt, mark=square]
        table{burg/epi6/med_cwyne1s.txt};
        \addlegendentry{CWY};
        \addplot[color=green!50!black, line width = 1.25pt, mark=square]
        table{burg/epi6/med_kiops.txt};
        \addlegendentry{ioCGS};
        \addplot[color=red!50!black, line width = 1.25pt, mark=square]table{burg/epi6/med_pmexne1s.txt};
        \addlegendentry{GSMGS};
        \addplot[color=yellow!50!black, line width = 1.25pt, mark=square]table{burg/epi6/med_icwyne1s.txt};
        \addlegendentry{NCWY};
        \end{axis}
        \end{tikzpicture}
          \caption{Strong scaling results for Burger's with low-synchronization hybrid methods for \texttt{epi6}.}
          \label{fig:burg_e6}
        \end{figure}
        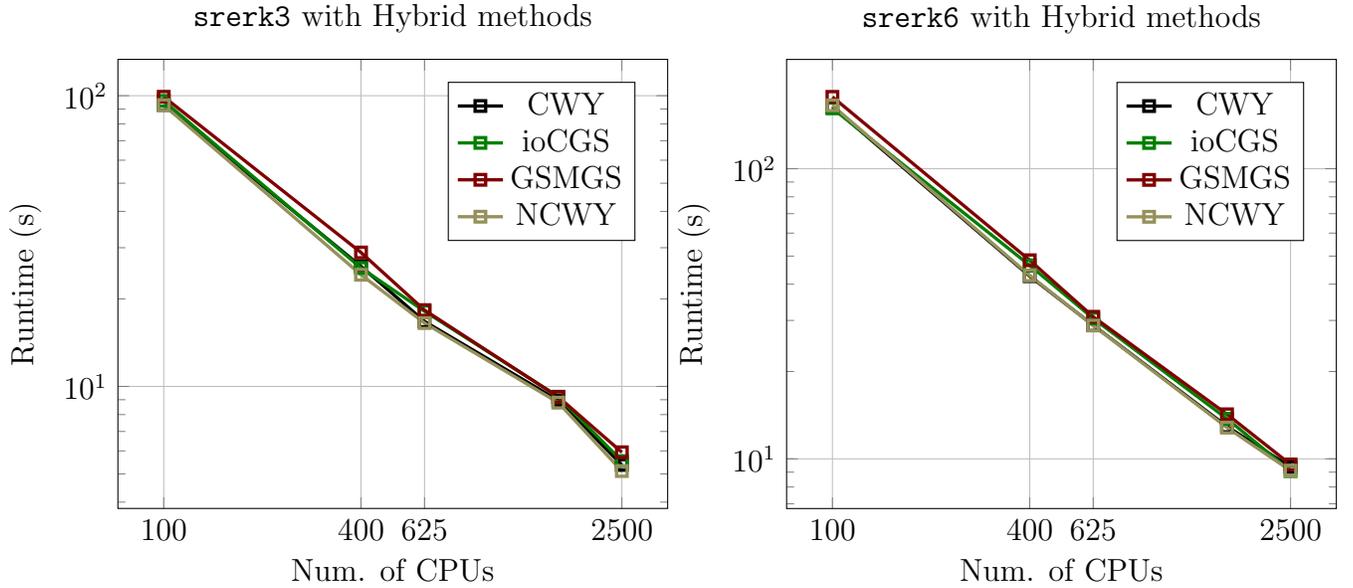
\begin{figure}[H]
          \centering
          \begin{subfigure}[b]{0.5\textwidth}
           \begin{tikzpicture}
        \begin{axis}[
        grid=major, width=1\textwidth, height=0.85\textwidth,
        xtick={100, 400, 625, 2500},
        xticklabels={100, 400, 625, 2500},
        xlabel={Num. of CPUs}, ylabel={Runtime (s)},
        title={\texttt{srerk3} with Hybrid methods},
        xmode=log, ymode=log,
        legend style={at={(0.60,0.95)},anchor=north west}
        ]
        \addplot[color=black!80!black, line width = 1.25pt, mark=square]
        table{burg/srerk3/med_cwyne1s.txt};
        \addlegendentry{CWY};
        \addplot[color=green!50!black, line width = 1.25pt, mark=square]
        table{burg/srerk3/med_kiops.txt};
        \addlegendentry{ioCGS};
        \addplot[color=red!50!black, line width = 1.25pt, mark=square]table{burg/srerk3/med_pmexne1s.txt};
        \addlegendentry{GSMGS};
        \addplot[color=yellow!50!black, line width = 1.25pt, mark=square]table{burg/srerk3/med_icwyne1s.txt};
        \addlegendentry{NCWY};
        \end{axis}
        \end{tikzpicture}
          \end{subfigure}%
         \hfill
          \begin{subfigure}[b]{0.5\textwidth}
              \begin{tikzpicture}
        \begin{axis}[
        grid=major, width=1\textwidth, height=0.85\textwidth,
        xtick={100, 400, 625, 2500},
        xticklabels={100, 400, 625, 2500},
        xlabel={Num. of CPUs}, ylabel={Runtime (s)},
        title={\texttt{srerk6} with Hybrid methods},
        xmode=log, ymode=log,
        legend style={at={(0.60,0.95)},anchor=north west}
        ]
        \addplot[color=black!80!black, line width = 1.25pt, mark=square]
        table{burg/srerk6/med_cwyne1s.txt};
        \addlegendentry{CWY};
        \addplot[color=green!50!black, line width = 1.25pt, mark=square]
        table{burg/srerk6/med_kiops.txt};
        \addlegendentry{ioCGS};
        \addplot[color=red!50!black, line width = 1.25pt, mark=square]table{burg/srerk6/med_pmexne1s.txt};
        \addlegendentry{GSMGS};
        \addplot[color=yellow!50!black, line width = 1.25pt, mark=square]table{burg/srerk6/med_icwyne1s.txt};
        \addlegendentry{NCWY};
        \end{axis}
        \end{tikzpicture}
          \end{subfigure}
          \caption{Strong scaling results for Burger's with low-synchronization hybrid methods for \texttt{srerk} integrators.}
          \label{fig:burg_srerk}
        \end{figure}

\begin{table}[H]
	\centering
	\begin{tabular}{c|c|c|c|c|c|c p{1.5cm}p{1.5cm}p{1.5cm}p{1.5cm}p{1.5cm}p{1cm}p{1cm}|}
		\toprule
		low-synchronization Method & $p$ & EPI4 & EPI5 & EPI6 & SRERK3 & SRERK6  \\ \hline
		\multirow{ 2}{*}{\texttt{H-GSMGS}} 
		& 100 & 0.98 & 0.98 & 1.00 & 1.03 & 1.10  \\
		& 400 & 1.01 & 0.95  & 1.00 & 1.13 & 1.04  \\
		& 625 & 1.09 & 0.95 & 1.11 & 1.01 & 1.01\\
		& 1600 & 1.13 & 0.96 & 0.99 & 1.00 & 1.04\\
		& 2500 & 1.01 & 1.02 & 1.14 & 1.07 & 1.05 \\
		\midrule
		\multirow{ 6}{*}{\texttt{H-NCWY}} 
		& 100 & 0.97 & 0.98 & 1.00 & 0.97 & 1.02 \\
		& 400 & 1.00 & 0.93  & 0.94 & 0.95 & 0.92 \\ 
		& 625 & 1.02 & 0.92 & 1.06 &  0.91 & 0.95 \\
		& 1600 & 1.18 & 1.01 & 0.95 & 0.95 & 0.93 \\
		& 2500 & 0.77 & 1.05 & 0.98 & 0.92 & 1.00 \\
		\midrule
		\multirow{ 6}{*}{\texttt{H-CWY}} 
		& 100 & 0.98 & 0.99 &  0.99 & 0.97 & 1.01  \\
		& 400 & 1.00 & 0.97  & 0.95 & 1.01 & 0.92 \\
		& 625 & 0.99 & 0.95 & 1.07 & 0.92 & 0.95\\
		& 1600 & 1.17 & 0.88 & 0.98 & 0.97 & 0.94 \\
		& 2500 & 0.75 & 1.03 & 1.07 & 0.97 & 1.03 \\
		\midrule 
	\end{tabular}
	\caption{Runtime ratios of \texttt{low-synchronization}/\texttt{ioCGS} for invisicid Burger's pde.}
	\label{runtimetratioburg}
\end{table}


    \subsection{Large Scale Shallow-Water Equations}\label{sec:SW}

        We now  present the results for the shallow water (SW) equations on the cubed-sphere, as presented in \cite{val_cube}. The SW equations are frequently used as a benchmark for assessing numerical methods in geophysical fluid dynamics and serve as a simplified model of choice for numerical weather prediction \cite{swbenchmark}. 
        
        We solve the SW equations using a direct flux reconstruction (DFR) spatial discretization \cite{val_cube,huynh2020discontinuous, romero2017development}. The computational grid is made up of $6 \times N_e \times N_e$ elements, each element with $N_s \times N_s$ solution points. We use $N_e = 168$ elements with $N_s = 7$ grid points, for a total of approximately 8.2 million degrees of freedom (8 ,297 ,856 total grid points). All simulations are carried out with a timestep of 900 seconds, a Krylov tolerance within \texttt{KIOPS} of $10^{-12}$, and a final time of 4 hours. The strong scaling results are carried out for processor sizes of $p = [ 864, 1176, 2646, 4704, 10584]$\footnote{For brevity, when referencing a specific processor size, the short hand $k$ (short for kilo which denotes 1000) will be used for size, i.e. $2k$ instead of $2646$} on the following standard benchmark cases:
        
        \begin{itemize}
             \item \textbf{SW-C5}: Zonal flow over isolated mountain (case 5 from Williamson et. al. \cite{swbenchmark})
            \item \textbf{SW-C6}: Rossby-Haurwitz wavenumber-4 (case 6 from Williamson et. al. \cite{swbenchmark}) 
            \item \textbf{SW-C8}: Unstable jet-stream from Galewsky et. al. \cite{galewsky2004initial}
            
        \end{itemize}

        \subsubsection{Strong Scaling Results}

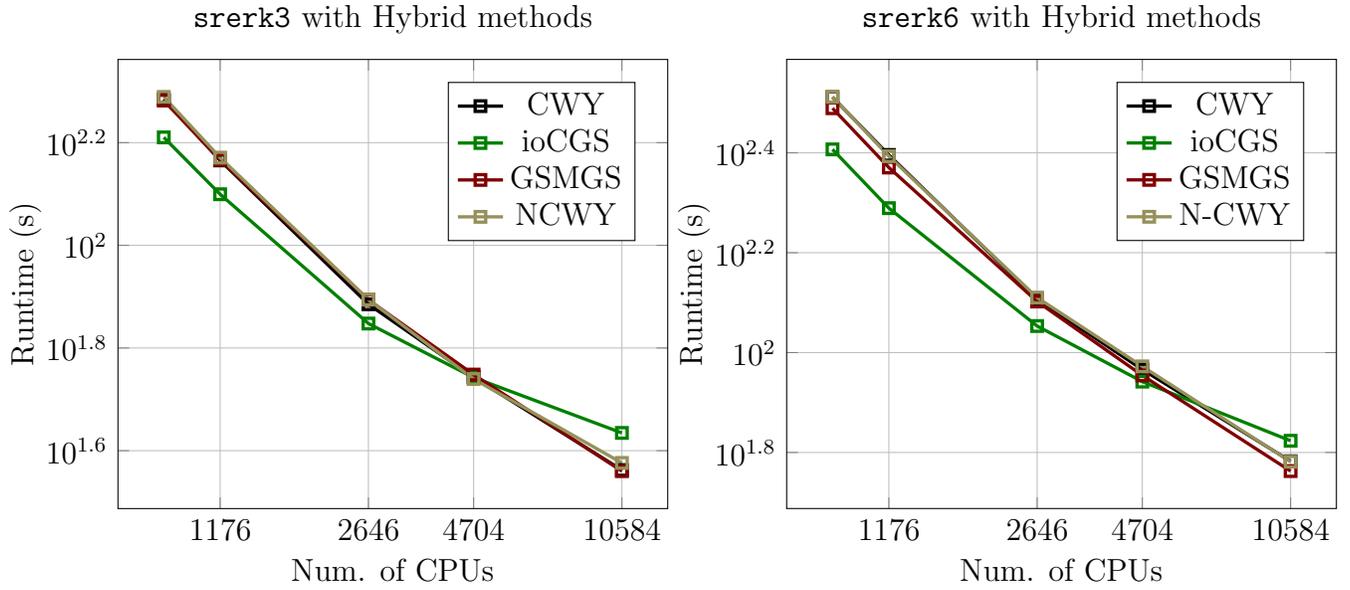
\begin{figure}[H]
          \centering
          \begin{subfigure}[b]{0.5\textwidth}
           \begin{tikzpicture}
        \begin{axis}[
        grid=major, width=1\textwidth, height=0.85\textwidth,
        xtick={ 1176, 2646, 4704, 10584},
        xticklabels={ 1176, 2646, 4704, 10584},
        xlabel={Num. of CPUs}, ylabel={Runtime (s)},
        title={\texttt{srerk3} with Hybrid methods},
        xmode=log, ymode=log,
        legend style={at={(0.60,0.95)},anchor=north west}
        ]
        \addplot[color=black!80!black, line width = 1.25pt, mark=square]
        table{srerk3/med_cwyne1s_c5.txt};
        \addlegendentry{CWY};
        \addplot[color=green!50!black, line width = 1.25pt, mark=square]
        table{srerk3/med_kiops_c5.txt};
        \addlegendentry{ioCGS};
        \addplot[color=red!50!black, line width = 1.25pt, mark=square]table{srerk3/med_pmexne1s_c5.txt};
        \addlegendentry{GSMGS};
        \addplot[color=yellow!50!black, line width = 1.25pt, mark=square]table{srerk3/med_icwyne1s_c5.txt};
        \addlegendentry{NCWY};
        \end{axis}
        \end{tikzpicture}
          \end{subfigure}%
         \hfill
          \begin{subfigure}[b]{0.5\textwidth}
              \begin{tikzpicture}
        \begin{axis}[
        grid=major, width=1\textwidth, height=0.85\textwidth,
        xtick={ 1176, 2646, 4704, 10584},
        xticklabels={ 1176, 2646, 4704, 10584},
        xlabel={Num. of CPUs}, ylabel={Runtime (s)},
        title={\texttt{srerk6} with Hybrid methods},
        xmode=log, ymode=log,
        legend style={at={(0.60,0.95)},anchor=north west}
        ]
        \addplot[color=black!80!black, line width = 1.25pt, mark=square]
        table{srerk6/med_cwyne1s_c5.txt};
        \addlegendentry{CWY};
        \addplot[color=green!50!black, line width = 1.25pt, mark=square]
        table{srerk6/med_kiops_c5.txt};
        \addlegendentry{ioCGS};
        \addplot[color=red!50!black, line width = 1.25pt, mark=square]table{srerk6/med_pmexne1s_c5.txt};
        \addlegendentry{GSMGS};
        \addplot[color=yellow!50!black, line width = 1.25pt, mark=square]table{srerk6/med_icwyne1s_c5.txt};
        \addlegendentry{N-CWY};
        \end{axis}
        \end{tikzpicture}
          \end{subfigure}
          \caption{Strong scaling results for \texttt{srerk3} and \texttt{srerk6} with hybrid low-synchronization methods for SW-C5.}
          \label{fig:srerk_c5}
        \end{figure}

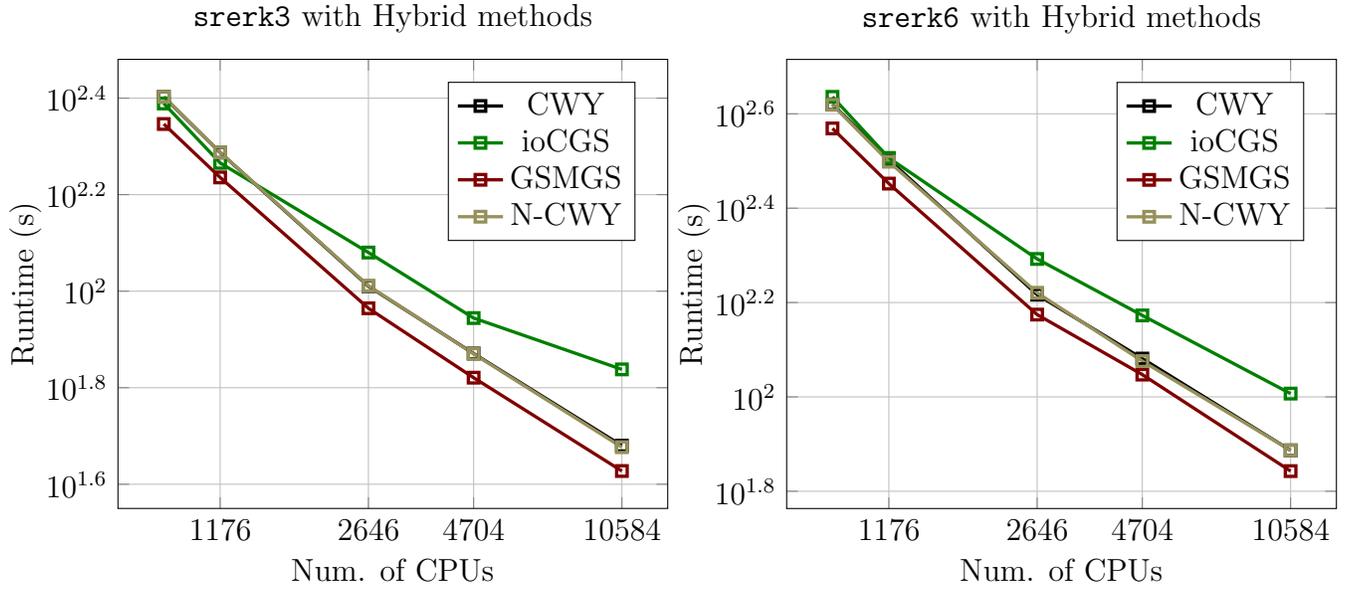
\begin{figure}[H]
          \centering
          \begin{subfigure}[b]{0.5\textwidth}
           \begin{tikzpicture}
        \begin{axis}[
        grid=major, width=1\textwidth, height=0.85\textwidth,
        xtick={ 1176, 2646, 4704, 10584},
        xticklabels={ 1176, 2646, 4704, 10584},
        xlabel={Num. of CPUs}, ylabel={Runtime (s)},
        title={\texttt{srerk3} with Hybrid methods},
        xmode=log, ymode=log,
        legend style={at={(0.60,0.95)},anchor=north west}
        ]
        \addplot[color=black!80!black, line width = 1.25pt, mark=square]
        table{srerk3/med_cwyne1s_c6.txt};
        \addlegendentry{CWY};
        \addplot[color=green!50!black, line width = 1.25pt, mark=square]
        table{srerk3/med_kiops_c6.txt};
        \addlegendentry{ioCGS};
        \addplot[color=red!50!black, line width = 1.25pt, mark=square]table{srerk3/med_pmexne1s_c6.txt};
        \addlegendentry{GSMGS};
        \addplot[color=yellow!50!black, line width = 1.25pt, mark=square]table{srerk3/med_icwyne1s_c6.txt};
        \addlegendentry{N-CWY};
        \end{axis}
        \end{tikzpicture}
          \end{subfigure}%
         \hfill
          \begin{subfigure}[b]{0.5\textwidth}
              \begin{tikzpicture}
        \begin{axis}[
        grid=major, width=1\textwidth, height=0.85\textwidth,
        xtick={ 1176, 2646, 4704, 10584},
        xticklabels={ 1176, 2646, 4704, 10584},
        xlabel={Num. of CPUs}, ylabel={Runtime (s)},
        title={\texttt{srerk6} with Hybrid methods},
        xmode=log, ymode=log,
        legend style={at={(0.60,0.95)},anchor=north west}
        ]
        \addplot[color=black!80!black, line width = 1.25pt, mark=square]
        table{srerk6/med_cwyne1s_c6.txt};
        \addlegendentry{CWY};
        \addplot[color=green!50!black, line width = 1.25pt, mark=square]
        table{srerk6/med_kiops_c6.txt};
        \addlegendentry{ioCGS};
        \addplot[color=red!50!black, line width = 1.25pt, mark=square]table{srerk6/med_pmexne1s_c6.txt};
        \addlegendentry{GSMGS};
        \addplot[color=yellow!50!black, line width = 1.25pt, mark=square]table{srerk6/med_icwyne1s_c6.txt};
        \addlegendentry{N-CWY};
        \end{axis}
        \end{tikzpicture}
          \end{subfigure}
          \caption{Strong scaling results for \texttt{srek3} and \texttt{srerk6} with hybrid low-synchronization methods for SW-C6.}
          \label{fig:srerk_c6}
        \end{figure}

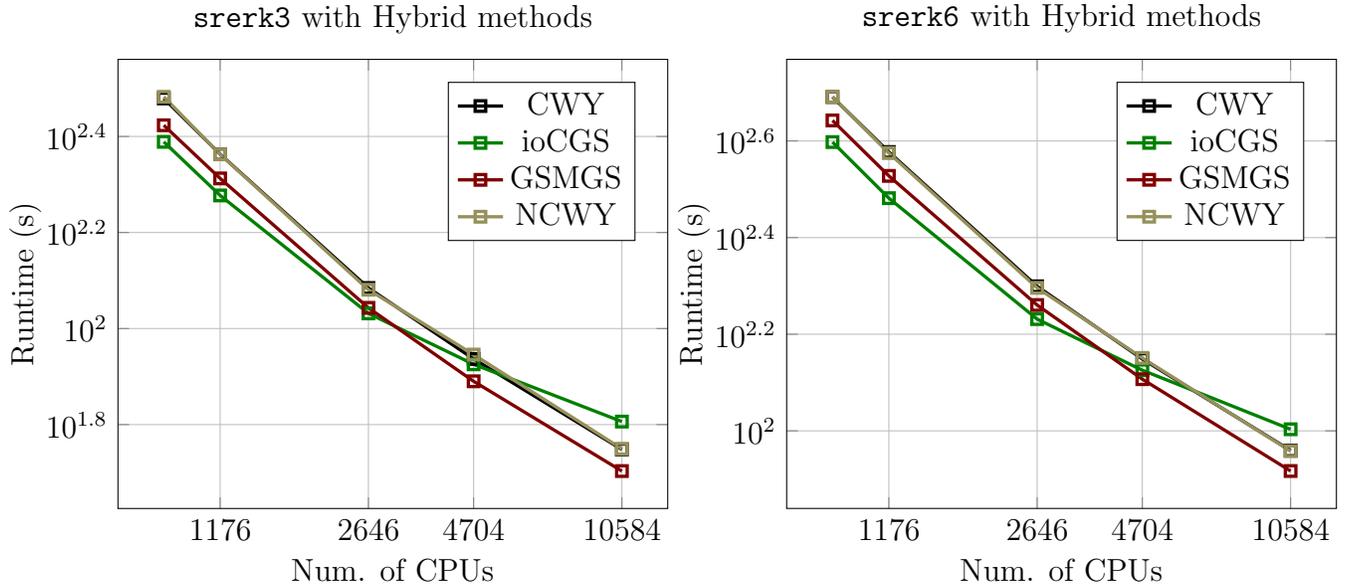
\begin{figure}[H]
          \centering
          \begin{subfigure}[b]{0.5\textwidth}
           \begin{tikzpicture}
        \begin{axis}[
        grid=major, width=1\textwidth, height=0.85\textwidth,
        xtick={ 1176, 2646, 4704, 10584},
        xticklabels={ 1176, 2646, 4704, 10584},
        xlabel={Num. of CPUs}, ylabel={Runtime (s)},
        title={\texttt{srerk3} with Hybrid methods},
        xmode=log, ymode=log,
        legend style={at={(0.60,0.95)},anchor=north west}
        ]
        \addplot[color=black!80!black, line width = 1.25pt, mark=square]
        table{srerk3/med_cwyne1s_c8.txt};
        \addlegendentry{CWY};
        \addplot[color=green!50!black, line width = 1.25pt, mark=square]
        table{srerk3/med_kiops_c8.txt};
        \addlegendentry{ioCGS};
        \addplot[color=red!50!black, line width = 1.25pt, mark=square]table{srerk3/med_pmexne1s_c8.txt};
        \addlegendentry{GSMGS};
        \addplot[color=yellow!50!black, line width = 1.25pt, mark=square]table{srerk3/med_icwyne1s_c8.txt};
        \addlegendentry{NCWY};
        \end{axis}
        \end{tikzpicture}
          \end{subfigure}%
         \hfill
          \begin{subfigure}[b]{0.5\textwidth}
              \begin{tikzpicture}
        \begin{axis}[
        grid=major, width=1\textwidth, height=0.85\textwidth,
        xtick={ 1176, 2646, 4704, 10584},
        xticklabels={ 1176, 2646, 4704, 10584},
        xlabel={Num. of CPUs}, ylabel={Runtime (s)},
        title={\texttt{srerk6} with Hybrid methods},
        xmode=log, ymode=log,
        legend style={at={(0.60,0.95)},anchor=north west}
        ]
        \addplot[color=black!80!black, line width = 1.25pt, mark=square]
        table{srerk6/med_cwyne1s_c8.txt};
        \addlegendentry{CWY};
        \addplot[color=green!50!black, line width = 1.25pt, mark=square]
        table{srerk6/med_kiops_c8.txt};
        \addlegendentry{ioCGS};
        \addplot[color=red!50!black, line width = 1.25pt, mark=square]table{srerk6/med_pmexne1s_c8.txt};
        \addlegendentry{GSMGS};
        \addplot[color=yellow!50!black, line width = 1.25pt, mark=square]table{srerk6/med_icwyne1s_c8.txt};
        \addlegendentry{NCWY};
        \end{axis}
        \end{tikzpicture}
          \end{subfigure}
          \caption{Strong scaling results for \texttt{srerk3} and \texttt{srerk6} with hybrid low-synchronization methods for SW-C8.}
          \label{fig:srerk_c8}
        \end{figure}

        The strong scaling results for SW-C5 are shown in Figures \ref{fig:epi45_c5} - \ref{fig:srerk_c5}, SW-C6 is shown in \ref{fig:epi45_c6} - \ref{fig:srerk_c6}, and SW-C8 is shown in \ref{fig:epi45_c8} - \ref{fig:srerk_c8} with the corresponding Tables \ref{runtimeratioc5}, \ref{runtimeratioc6}, \ref{runtimeratioc8} for SW-C5, SW-C6, and SW-C8 respectively.

        For these larger scale problems, we see the low-synchronization methods outperformed \texttt{ioCGS} for all cases when the number of processors was large enough. This is precisely what is expected given that for large number of processors the savings in communication that the low-synchronization methods offer play the key role in the overall time-to-solution of the exponential methods. The numerical test results support our hypothesis and highlight the advantages of the low-synchronization methods for very large scale stiff problems that have to be solved on massively parallel computing platforms. For this example set, the \texttt{H-GSMGS} algorithm showed the most improvement in performance. This could be explained by \texttt{H-GSMGS} having the largest computational complexity for computing the correction matrix $T$ (matrix-matrix product of $M^{-1}$ and $N$, and a lower triangular solve), as it is advantageous to increase the computational intensity in a parallel setting \cite{lund2023adaptively}.

        Of the three geophysical scenarios presented, exponential integrators for the Rossby-Haurwitz waves test problem (SW-C6) showed the most improvement in time-to-solution and scaling properties; see Figures \ref{fig:epi45_c6} - \ref{fig:srerk_c6} and Table \ref{runtimeratioc6}. First, all integrators showed an improved performance with the low-synchronization methods compared to \texttt{ioCGS} starting from the smallest processor size of $864$ CPUs. The greatest improvement over all integrators is a $45\%$ reduction in runtime for \texttt{epi6} with \texttt{H-GSMGS} for $10k$ processors. The integrators were consistently faster with \texttt{H-GSMGS} compared to \texttt{H-CWY} and \texttt{H-NCWY}, where the \texttt{H-WY} methods slowed down for $1k$ processors. Second, the integrators displayed better strong scaling properties with the low-synchronization methods. The \texttt{epi} methods display zig-zaggy behavior with \texttt{ioCGS} and show little improvement in runtime between $4k$ and $10k$ processors, while the strong scaling curves with the low-synchronization methods remains fairly linear. 
        
        The results for zonal flow over mountain (SW-C5) showed smaller gains but improvements were still observed once the number of processors exceeded $4k$; see Figures \ref{fig:epi45_c5} - \ref{fig:srerk_c5} and Table \ref{runtimeratioc5}. Exponential integrators with \texttt{ioCGS} remained more efficient for up to $4k$ processors for both \texttt{epi} and \texttt{srerk} type. For the \texttt{epi} integrators, the greatest improvement was seen with \texttt{epi5} with a $15\%$ reduction in runtime with \texttt{H-NCWY} and \texttt{H-CWY}. For the \texttt{srerk} integrators, the \texttt{srerk3} saw the best improvement with a $16\%$ reduction in runtime with \texttt{H-GSMGS}.
        
        The unstable jet stream (SW-C8) case presented intermediate improvements with the low-synchronization methods with all exponential integrators; see Figures \ref{fig:epi45_c8} - \ref{fig:srerk_c8} and Table \ref{runtimeratioc8}. The strong scaling graphs for the \texttt{epi} methods show a similar runtime between \texttt{ioCGS} and \texttt{H-GSMGS} for processors fewer than $2k$, followed by an improvement in performance by \texttt{H-GSMGS} as the number of processors increased (see Figures \ref{fig:epi45_c8} - \ref{fig:c8_e6}). Performance with the \texttt{H-WY} methods improved for more than $4k$ processors with \texttt{epi4} and the \texttt{srerk} type, and more than $2k$ for \texttt{epi5} and \texttt{epi6}. The largest improvement was for \texttt{epi6} with a $30\%$ reduction in runtime using \texttt{H-GSMGS} with $10k$ processors.  

        \begin{figure}[H]
          \centering
          \begin{subfigure}[b]{0.5\textwidth}
          \begin{tikzpicture}
        \begin{axis}[
        grid=major, width=1\textwidth, height=0.85\textwidth,
        title={\texttt{EPI4} with Hybrid methods},
        xtick={ 1176, 2646, 4704, 10584},
        xticklabels={ 1176, 2646, 4704, 10584},
        xlabel={Num. of CPUs}, ylabel={Runtime (s)},
        xmode=log, ymode=log,
        legend style={at={(0.60,0.95)},anchor=north west}
        ]
        \addplot[color=black!80!black, line width = 1.25pt, mark=square]
        table{epi4/med_cwyne1s_c5.txt};
        \addlegendentry{CWY};
        \addplot[color=green!50!black, line width = 1.25pt, mark=square]
        table{epi4/med_kiops_c5.txt};
        \addlegendentry{ioCGS};
        \addplot[color=red!50!black, line width = 1.25pt, mark=square]table{epi4/med_pmexne1s_c5.txt};
        \addlegendentry{GSMGS};
        \addplot[color=yellow!50!black, line width = 1.25pt, mark=square]table{epi4/med_icwyne1s_c5.txt};
        \addlegendentry{NCWY};
        \end{axis}
        \end{tikzpicture}
          \end{subfigure}%
         \hfill
          \begin{subfigure}[b]{0.5\textwidth}
              \begin{tikzpicture}
        \begin{axis}[
        grid=major, width=1\textwidth, height=0.85\textwidth,
        xtick={ 1176, 2646, 4704, 10584},
        xticklabels={ 1176, 2646, 4704, 10584},
        xlabel={Num. of CPUs}, ylabel={Runtime (s)},
        title={\texttt{EPI5} with Hybrid methods},
        xmode=log, ymode=log,
        legend style={at={(0.60,0.95)},anchor=north west}
        ]
        \addplot[color=black!80!black, line width = 1.25pt, mark=square]
        table{epi5/med_cwyne1s_c5.txt};
        \addlegendentry{CWY};
        \addplot[color=green!50!black, line width = 1.25pt, mark=square]
        table{epi5/med_kiops_c5.txt};
        \addlegendentry{ioCGS};
        \addplot[color=red!50!black, line width = 1.25pt, mark=square]table{epi5/med_pmexne1s_c5.txt};
        \addlegendentry{GSMGS};
        \addplot[color=yellow!50!black, line width = 1.25pt, mark=square]table{epi5/med_icwyne1s_c5.txt};
        \addlegendentry{NCWY};
        \end{axis}
        \end{tikzpicture}
          \end{subfigure}
          \caption{Strong scaling results for \texttt{epi4} and \texttt{epi5} with hybrid low-synchronization methods for SW-C5.}
          \label{fig:epi45_c5}
        \end{figure}
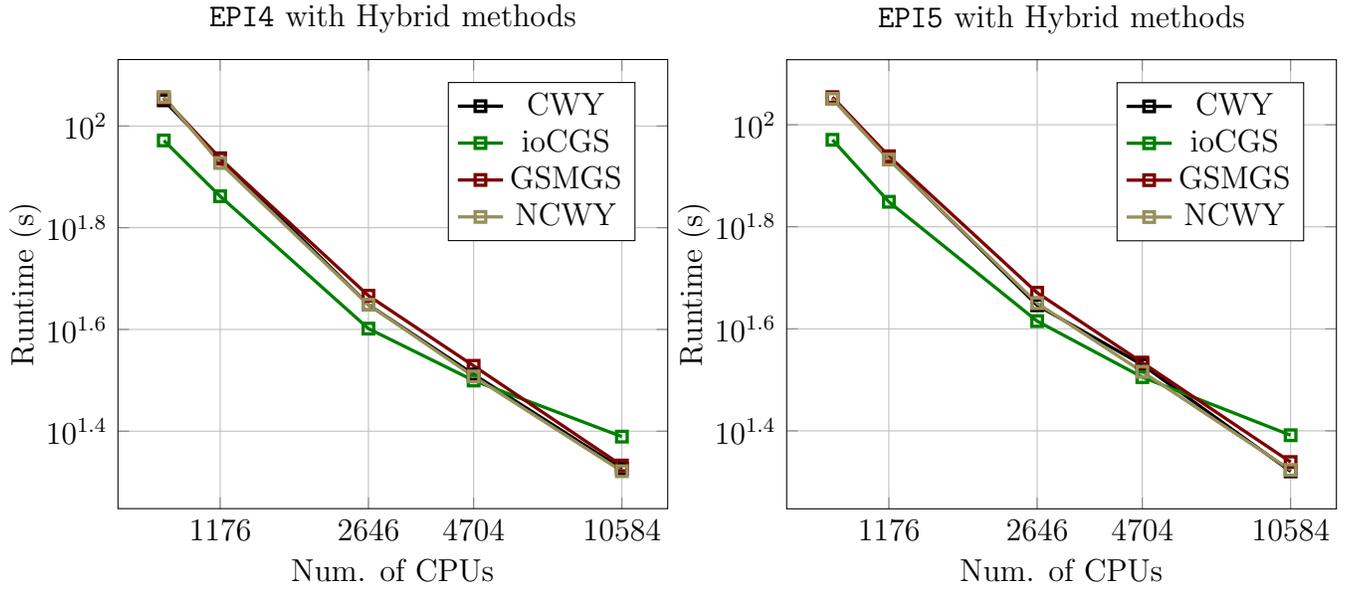
         \begin{figure}[H]
          \centering
          \begin{tikzpicture}
        \begin{axis}[
        grid=major, width=0.52\textwidth, height=0.43\textwidth,
        title={\texttt{EPI6} with Hybrid methods},
        xtick={1176, 2646, 4704, 10584},
        xticklabels={1176, 2646, 4704, 10584},
        xlabel={Num. of CPUs}, ylabel={Runtime (s)},
        xmode=log, ymode=log,
        legend style={at={(0.60,0.95)},anchor=north west}
        ]
        \addplot[color=black!80!black, line width = 1.25pt, mark=square]
        table{epi6/med_cwyne1s_c5.txt};
        \addlegendentry{CWY};
        \addplot[color=green!50!black, line width = 1.25pt, mark=square]
        table{epi6/med_kiops_c5.txt};
        \addlegendentry{ioCGS};
        \addplot[color=red!50!black, line width = 1.25pt, mark=square]table{epi6/med_pmexne1s_c5.txt};
        \addlegendentry{GSMGS};
        \addplot[color=yellow!50!black, line width = 1.25pt, mark=square]table{epi6/med_icwyne1s_c5.txt};
        \addlegendentry{NCWY};
        \end{axis}
        \end{tikzpicture}
          \caption{Strong scaling results for \texttt{epi6} with hybrid low-synchronization methods for SW-C5.}
          \label{fig:c5_e6}
        \end{figure}
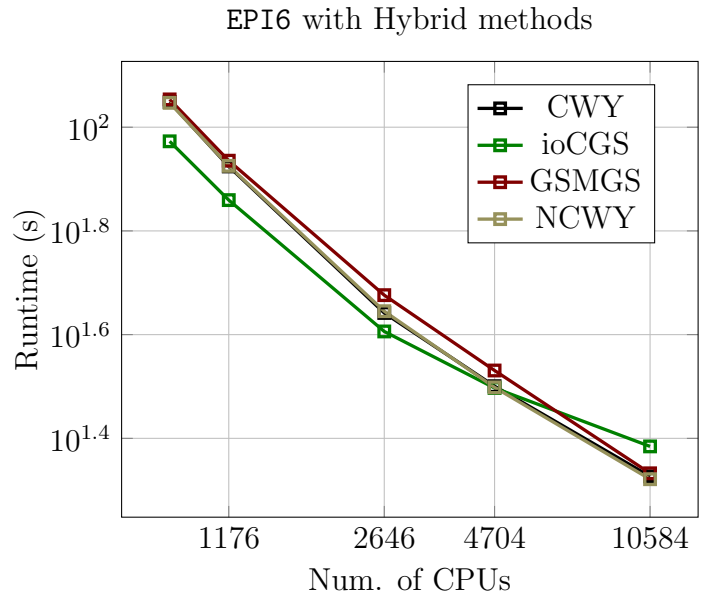

        \begin{table}[H]
	\centering
	\begin{tabular}{c|c|c|c|c|c|c p{1cm}p{1cm}p{1cm}p{1cm}p{1cm}p{1cm}p{1cm}|}
		\toprule
		low-synchronization Method & $p$ & EPI4 & EPI5 & EPI6 & SRERK3 & SRERK6  \\ \hline
		\multirow{ 2}{*}{\texttt{H-GSMGS}} 
		& 864 & 1.21 & 1.21  & 1.20 & 1.18 & 1.21 \\
		& 1k & 1.18 & 1.23 & 1.19 & 1.17 & 1.21\\
		& 2k & 1.16 & 1.14 & 1.17 & 1.11 & 1.12 \\
		& 4k & 1.07 & 1.07 & 1.08 & 1.01 & 1.03 \\
		& 10k & 0.88 & 0.89 & 0.89 & 0.84 & 0.87 \\
		\midrule
		\multirow{ 6}{*}{\texttt{H-NCWY}} 
		& 864 & 1.22 & 1.20  & 1.19 & 1.20 &  1.27 \\ 
		& 1k & 1.16 & 1.21 & 1.16 & 1.18 & 1.27\\
		& 2k & 1.11 & 1.08 & 1.09 & 1.11 & 1.14 \\
		& 4k & 1.02 & 1.02 & 1.00 & 0.99 & 1.07 \\
		& 10k & 0.86 & 0.85 & 0.86 & 0.87 & 0.91\\
		\midrule
		\multirow{ 6}{*}{\texttt{H-CWY}} 
		& 864 & 1.20 & 1.20  & 1.19 & 1.19 & 1.27 \\
		& 1k & 1.19 & 1.21 & 1.16 & 1.16 & 1.28\\
		& 2k & 1.11 & 1.07 & 1.09 & 1.09 & 1.13\\
		& 4k & 1.03 & 1.06 & 1.01 & 1.00 & 1.06\\
		& 10k & 0.86 & 0.85 & 0.88 & 0.85 & 0.91 \\
		\midrule 
	\end{tabular}
	\caption{Runtime ratios of \texttt{low-synchronization}/\texttt{ioCGS} for SW-C5. }
	\label{runtimeratioc5}
        \end{table}

        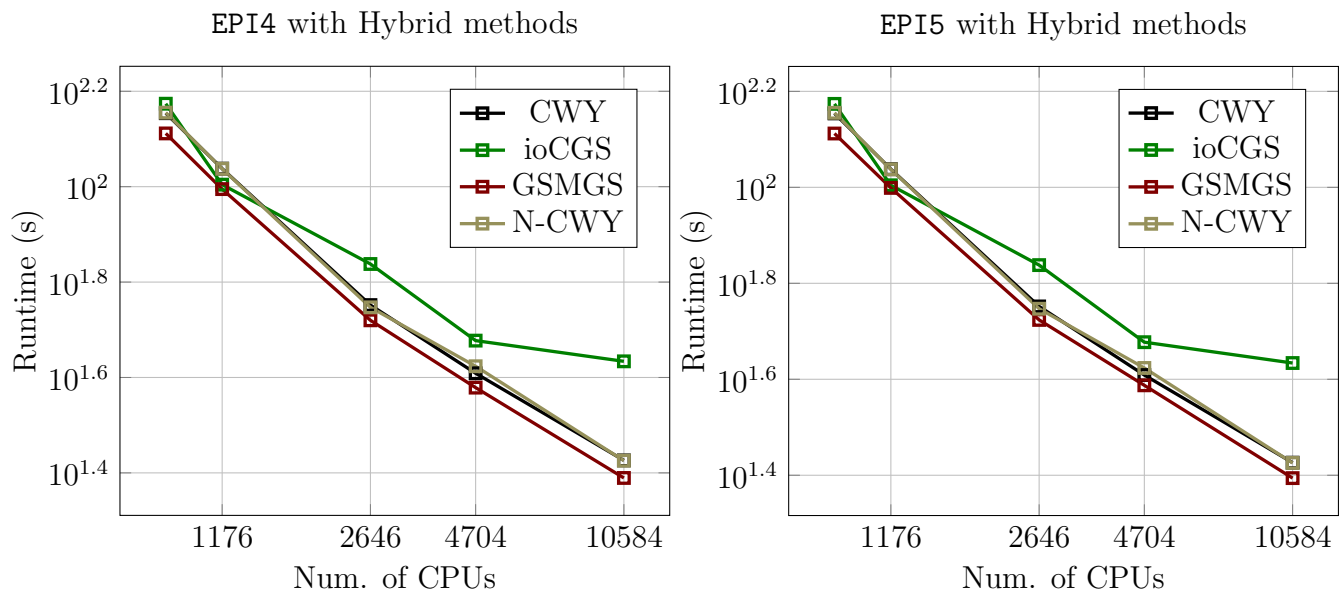
\begin{figure}[H]
          \centering
          \begin{subfigure}[b]{0.5\textwidth}
          \begin{tikzpicture}
        \begin{axis}[
        grid=major, width=1\textwidth, height=0.85\textwidth,
        title={\texttt{EPI4} with Hybrid methods},
        xtick={1176, 2646, 4704, 10584},
        xticklabels={1176, 2646, 4704, 10584},
        xlabel={Num. of CPUs}, ylabel={Runtime (s)},
        xmode=log, ymode=log,
        legend style={at={(0.60,0.95)},anchor=north west}
        ]
        \addplot[color=black!80!black, line width = 1.25pt, mark=square]
        table{epi4/med_cwyne1s_c6.txt};
        \addlegendentry{CWY};
        \addplot[color=green!50!black, line width = 1.25pt, mark=square]
        table{epi4/med_kiops_c6.txt};
        \addlegendentry{ioCGS};
        \addplot[color=red!50!black, line width = 1.25pt, mark=square]table{epi4/med_pmexne1s_c6.txt};
        \addlegendentry{GSMGS};
        \addplot[color=yellow!50!black, line width = 1.25pt, mark=square]table{epi4/med_icwyne1s_c6.txt};
        \addlegendentry{N-CWY};
        \end{axis}
        \end{tikzpicture}
          \end{subfigure}%
         \hfill
          \begin{subfigure}[b]{0.5\textwidth}
              \begin{tikzpicture}
        \begin{axis}[
        grid=major, width=1\textwidth, height=0.85\textwidth,
        xtick={1176, 2646, 4704, 10584},
        xticklabels={1176, 2646, 4704, 10584},
        xlabel={Num. of CPUs}, ylabel={Runtime (s)},
        title={\texttt{EPI5} with Hybrid methods},
        xmode=log, ymode=log,
        legend style={at={(0.60,0.95)},anchor=north west}
        ]
        \addplot[color=black!80!black, line width = 1.25pt, mark=square]
        table{epi4/med_cwyne1s_c6.txt};
        \addlegendentry{CWY};
        \addplot[color=green!50!black, line width = 1.25pt, mark=square]
        table{epi4/med_kiops_c6.txt};
        \addlegendentry{ioCGS};
        \addplot[color=red!50!black, line width = 1.25pt, mark=square]table{epi5/med_pmexne1s_c6.txt};
        \addlegendentry{GSMGS};
        \addplot[color=yellow!50!black, line width = 1.25pt, mark=square]table{epi4/med_icwyne1s_c6.txt};
        \addlegendentry{N-CWY};
        \end{axis}
        \end{tikzpicture}
          \end{subfigure}
          \caption{Strong scaling results for \texttt{epi4} and \texttt{epi5} with hybrid low-synchronization methods for SW-C6.}
          \label{fig:epi45_c6}
        \end{figure}
        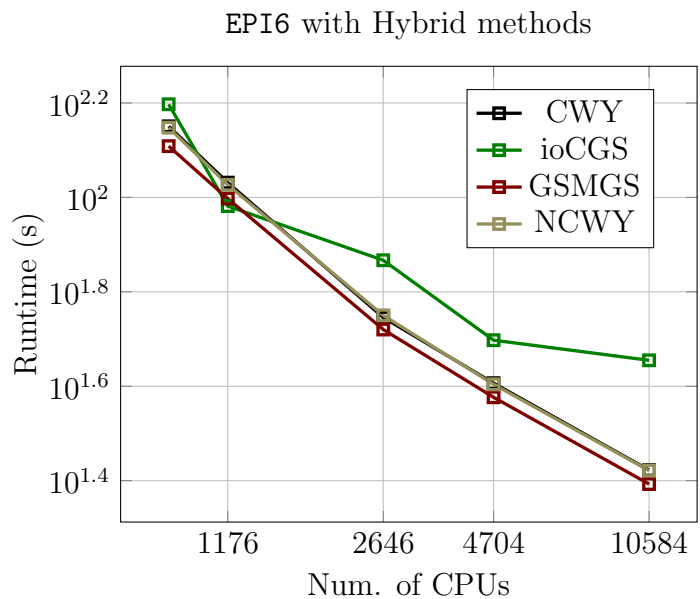
\begin{figure}[H]
          \centering
          \begin{tikzpicture}
        \begin{axis}[
        grid=major, width=0.52\textwidth, height=0.43\textwidth,
        title={\texttt{EPI6} with Hybrid methods},
        xtick={1176, 2646, 4704, 10584},
        xticklabels={1176, 2646, 4704, 10584},
        xlabel={Num. of CPUs}, ylabel={Runtime (s)},
        xmode=log, ymode=log,
        legend style={at={(0.60,0.95)},anchor=north west}
        ]
        \addplot[color=black!80!black, line width = 1.25pt, mark=square]
        table{epi6/med_cwyne1s_c6.txt};
        \addlegendentry{CWY};
        \addplot[color=green!50!black, line width = 1.25pt, mark=square]
        table{epi6/med_kiops_c6.txt};
        \addlegendentry{ioCGS};
        \addplot[color=red!50!black, line width = 1.25pt, mark=square]table{epi6/med_pmexne1s_c6.txt};
        \addlegendentry{GSMGS};
        \addplot[color=yellow!50!black, line width = 1.25pt, mark=square]table{epi6/med_icwyne1s_c6.txt};
        \addlegendentry{NCWY};
        \end{axis}
        \end{tikzpicture}
          \caption{Strong scaling results for \texttt{epi6} with hybrid low-synchronization methods for SW-C6.}
          \label{fig:c6_e6}
        \end{figure}

        \begin{table}[H]
	\centering
	\begin{tabular}{c|c|c|c|c|c|c p{1cm}p{1cm}p{1cm}p{1cm}p{1cm}p{1cm}p{1cm}|}
		\toprule
		low-synchronization Method & $p$ & EPI4 & EPI5 & EPI6 & SRERK3 & SRERK6  \\ \hline
		\multirow{ 2}{*}{\texttt{H-GSMGS}} 
		& 864 & 0.86 & 0.84  & 0.81 & 0.91 & 0.86 \\
		& 1k & 0.98 & 1.02 & 1.04 & 0.93 & 0.88\\
		& 2k & 0.76 & 0.74 & 0.71 & 0.77 & 0.76 \\
		& 4k & 0.78 & 0.81 & 0.76 & 0.75 & 0.75 \\
		& 10k & 0.57 & 0.56 & 0.55 & 0.62 & 0.69 \\
		\midrule
		\multirow{ 6}{*}{\texttt{H-NCWY}} 
		& 864 & 0.96 & 0.92  & 0.89 & 1.03 &  0.96 \\ 
		& 1k & 1.08 & 1.10 & 1.11 & 1.05 & 0.98\\
		& 2k & 0.81 & 0.80 & 0.76 & 0.85 & 0.85 \\
		& 4k & 0.88 & 0.84 & 0.81 & 0.84 & 0.80 \\
		& 10k & 0.62 & 0.60 & 0.58 & 0.69 & 0.76\\
		\midrule
		\multirow{ 6}{*}{\texttt{H-CWY}} 
		& 864 & 0.96 & 0.92  & 0.90 & 1.03 & 0.96\\
		& 1k & 1.08 & 1.11 & 1.12 & 1.05 & 0.99\\
		& 2k & 0.82 & 0.80 & 0.76 & 0.85 & 0.84\\
		& 4k & 0.85 & 0.85 & 0.81 & 0.84 & 0.81\\
		& 10k & 0.62 & 0.60 & 0.59 & 0.69 & 0.76 \\
		\midrule 
	\end{tabular}
	\caption{Runtime ratios of \texttt{low-synchronization}/\texttt{ioCGS} for SW-C6.}
	\label{runtimeratioc6}
        \end{table}

        \begin{figure}[H]
          \centering
          \begin{subfigure}[b]{0.5\textwidth}
          \begin{tikzpicture}
        \begin{axis}[
        grid=major, width=1\textwidth, height=0.85\textwidth,
        title={\texttt{EPI4} with Hybrid methods},
        xtick={ 1176, 2646, 4704, 10584},
        xticklabels={ 1176, 2646, 4704, 10584},
        xlabel={Num. of CPUs}, ylabel={Runtime (s)},
        xmode=log, ymode=log,
        legend style={at={(0.60,0.95)},anchor=north west}
        ]
        \addplot[color=black!80!black, line width = 1.25pt, mark=square]
        table{epi4/med_cwyne1s_c8.txt};
        \addlegendentry{CWY};
        \addplot[color=green!50!black, line width = 1.25pt, mark=square]
        table{epi4/med_kiops_c8.txt};
        \addlegendentry{ioCGS};
        \addplot[color=red!50!black, line width = 1.25pt, mark=square]table{epi4/med_pmexne1s_c8.txt};
        \addlegendentry{GSMGS};
        \addplot[color=yellow!50!black, line width = 1.25pt, mark=square]table{epi4/med_icwyne1s_c8.txt};
        \addlegendentry{NCWY};
        \end{axis}
        \end{tikzpicture}
          \end{subfigure}%
         \hfill
          \begin{subfigure}[b]{0.5\textwidth}
              \begin{tikzpicture}
        \begin{axis}[
        grid=major, width=1\textwidth, height=0.85\textwidth,
        xtick={ 1176, 2646, 4704, 10584},
        xticklabels={ 1176, 2646, 4704, 10584},
        xlabel={Num. of CPUs}, ylabel={Runtime (s)},
        title={\texttt{EPI5} with \texttt{NE}-methods},
        xmode=log, ymode=log,
        legend style={at={(0.60,0.95)},anchor=north west}
        ]
        \addplot[color=black!80!black, line width = 1.25pt, mark=square]
        table{epi5/med_cwyne1s_c8.txt};
        \addlegendentry{CWY};
        \addplot[color=green!50!black, line width = 1.25pt, mark=square]
        table{epi5/med_kiops_c8.txt};
        \addlegendentry{ioCGS};
        \addplot[color=red!50!black, line width = 1.25pt, mark=square]table{epi5/med_pmexne1s_c8.txt};
        \addlegendentry{GSMGS};
        \addplot[color=yellow!50!black, line width = 1.25pt, mark=square]table{epi5/med_icwyne1s_c8.txt};
        \addlegendentry{NCWY};
        \end{axis}
        \end{tikzpicture}
          \end{subfigure}
          \caption{Strong scaling results for \texttt{epi4} and \texttt{epi5} with hybrid low-synchronization methods for SW-C8.}
          \label{fig:epi45_c8}
        \end{figure}
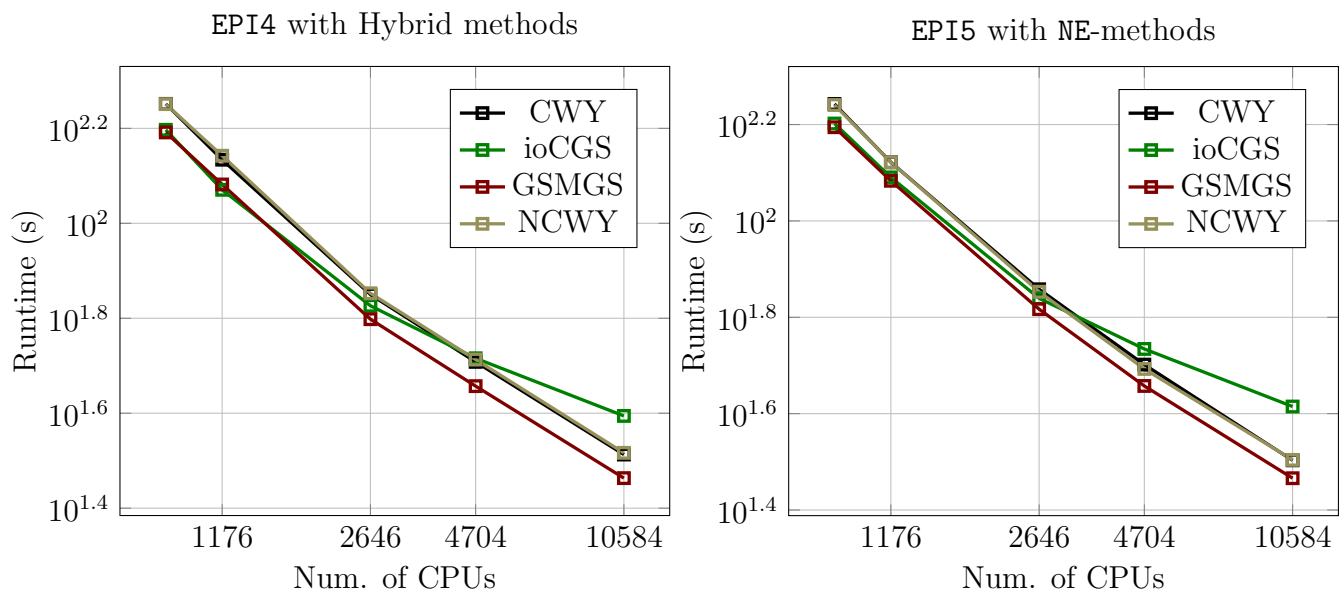
        \begin{figure}[H]
          \centering
          \begin{tikzpicture}
        \begin{axis}[
        grid=major, width=0.52\textwidth, height=0.43\textwidth,
        title={\texttt{EPI6} with Hybrid methods},
        xtick={ 1176, 2646, 4704, 10584},
        xticklabels={ 1176, 2646, 4704, 10584},
        xlabel={Num. of CPUs}, ylabel={Runtime (s)},
        xmode=log, ymode=log,
        legend style={at={(0.60,0.95)},anchor=north west}
        ]
        \addplot[color=black!80!black, line width = 1.25pt, mark=square]
        table{epi6/med_cwyne1s_c8.txt};
        \addlegendentry{CWY};
        \addplot[color=green!50!black, line width = 1.25pt, mark=square]
        table{epi6/med_kiops_c8.txt};
        \addlegendentry{ioCGS};
        \addplot[color=red!50!black, line width = 1.25pt, mark=square]table{epi6/med_pmexne1s_c8.txt};
        \addlegendentry{GSMGS};
        \addplot[color=yellow!50!black, line width = 1.25pt, mark=square]table{epi6/med_icwyne1s_c8.txt};
        \addlegendentry{NCWY};
        \end{axis}
        \end{tikzpicture}
          \caption{Strong scaling results for \texttt{epi6} with hybrid low-synchronization methods for SW-C8.}
          \label{fig:c8_e6}
        \end{figure}
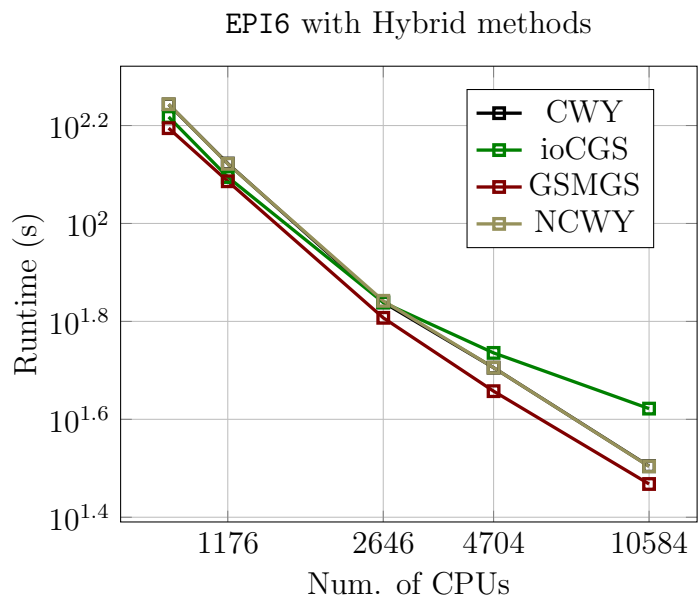
        
        \begin{table}[H]
	\centering
	\begin{tabular}{c|c|c|c|c|c|c p{1cm}p{1cm}p{1cm}p{1cm}p{1cm}p{1cm}p{1cm}|}
		\toprule
		low-synchronization Method & $p$ & EPI4 & EPI5 & EPI6 & SRERK3 & SRERK6  \\ \hline
		\multirow{ 2}{*}{\texttt{H-GSMGS}} 
		& 864 & 0.99 & 0.98  & 0.95 & 1.08 & 1.11 \\
		& 1k & 1.03 & 0.98 & 0.98 & 1.08 & 1.11\\
		& 2k & 0.94 & 0.95 & 0.93 & 1.03 & 1.07 \\
		& 4k & 0.87 & 0.84 & 0.84 & 0.92 & 0.96 \\
		& 10k & 0.74 & 0.71 & 0.70 & 0.79 & 0.82 \\
		\midrule
		\multirow{ 6}{*}{\texttt{H-NCWY}} 
		& 864 & 1.13 & 1.09  & 1.06 & 1.24 &  1.24 \\ 
		& 1k & 1.18 & 1.08 & 1.06 & 1.22 & 1.24\\
		& 2k & 1.06 & 1.03 & 1.01 & 1.12 & 1.16 \\
		& 4k & 0.99 & 0.91 & 0.93 & 1.05 & 1.06 \\
		& 10k & 0.84 & 0.77 & 0.76 & 0.88 & 0.90\\
		\midrule
		\multirow{ 6}{*}{\texttt{H-CWY}} 
		& 864 & 1.13 & 1.10  & 1.06 & 1.23 & 1.24 \\
		& 1k & 1.16 & 1.07 & 1.07 & 1.22 & 1.25\\
		& 2k & 1.06 & 1.04 & 1.00 & 1.13 & 1.17\\
		& 4k & 0.98 & 0.93 & 0.93 & 1.03 & 1.05\\
		& 10k & 0.83 & 0.77 & 0.76 & 0.87 & 0.90 \\
		\midrule 
	\end{tabular}
	\caption{Runtime ratios of \texttt{low-synchronization}/\texttt{ioCGS} for SW-C8.}
	\label{runtimeratioc8}
        \end{table}


    \section{Conclusion}

    While Krylov projection-based algorithms are the methods of choice when a large scale stiff problem is integrated in time with an exponential integrator, the global synchronizations required to compute inner products at each Arnoldi iteration impose limitations on parallel scalability of such exponential-Krylov techniques. For example, using Modified Gram-Schmidt for the orthonalization of the Krylov basis requires a significant number of global synchronizations that is proportional to the size of the Krylov subspace. In \texttt{KIOPS}, the Krylov subspace has to be  recomputed at each substep $e^{\tau_i A}b_{\tau_i}$ (\ref{subeqn:substep} in sect. \ref{sec:KryApprox}). Therefore, each substep relies on the Arnoldi iteration as the core component that determines the computational efficiency of exponential integrators. While incomplete orthogonalization can ameliorate this limitations for parallel implementation of an exponential integrator, many problems require full orthogonalization to be performed for reasons of accuracy and efficiency \cite{stewart2023variable}. In this paper, we introduced full orthogonalization, hybrid low-synchronization MGS variants for the Krylov approximation of the matrix exponential for use inside an exponential integrators. These methods reduce the number of global synchronizations to one per Arnoldi-step independent of the Krylov dimension.
    
    We studied performance of the new hybrid low-synchronization methods using a set of six test problems including three standard examples for testing stiff integrators and three problems that model different geophysical scenarios with the shallow water equations on the cubed-sphere grid. The numerical experiments confirmed the effectiveness of the low-synchronization approach for very large scale systems solved on massively parallel computing platforms. Exponential integrators with low-synchronization algorithms showed improved strong scalability with large numbers of processors in a CPU-only parallel environment. While the degree of improvement varied with the test cases and the integrators, the fact that the performance was improved for the largest number of processors remained consistent. The low-synchronization methods are designed for massively parallel systems such as multi-GPU implementations \cite{wy, postmodern}, therefore it follows that the most improvement is seen with the largest processor size; $2k$ for sect. \ref{sec:expdes}  and $10k$ for shallow water (sect. \ref{sec:SW}). Based on the results from the example problem in section \ref{sec:expdes}, it is unclear which is the better low-synchronization algorithm to pair with exponential integrators. However, for the shallow water examples, specifically SW-C8 (fig. \ref{fig:epi45_c8} - \ref{fig:srerk_c8}) and SW-C6 (fig. \ref{fig:epi45_c6} -  \ref{fig:srerk_c6}), the hybrid Iterated Gau{\ss}--Seidel MGS algorithm (\texttt{H-GSMGS} ) (algorithm \ref{alg:lowsync_skeleton} and \ref{alg:Tgsmsg} ) had the greatest runtime improvements. 

    Although the numerical experiments compare an incomplete orthogonalization CGS and full orthogonalization MGS, the accuracy of each substep projection $e^{\tau_i A}b_{\tau_i}$ and the overall solution is not expected to improve. The adaptivity in the \texttt{KIOPS} algorithm for $m_i$ and $\tau_i$ will ensure that at each substep the Krylov approximation $e^{\tau_i A}b_{\tau_i}$ has converged to the accuracy of the user-specified tolerance. However, with full-orthogonalization schemes, it is possible to have better substepping properties such as requiring less substeps. Therefore, in future work we aim to continue to improve the efficiency of the Krylov approximation by shifting towards the substepping procedure mentioned in section \ref{sec:KryApprox}. For example, incorporating higher order adaptive time-stepping automatic controllers for $\tau$ such as proportional-integral (PI) or proportional-integral-derivative (PID) controllers based off of \cite{soderlind2002automatic, soderlind2003digital, gustafsson1988pi}. High-order adaptivity techniques for the substep size $\tau$ could lead to reducing the number of steps and/or the number of rejected steps. Refining the substepping procedure could lead to a further improvement in performance of exponential integrators.

    \section{Code availability}

        The code used to generate the results is available at: \textbf{create a repo with my branch of WX-factory}.

    \section{Acknowledgements}
    The origin of this study goes back to fruitful discussions with Stephen Thomas. Stimulating discussions with Kathryn Lund are also acknowledged. The authors are grateful to Vincent Magnoux for his assistance in running parallel simulations on the supercomputer which was crucial to the success of this study. This work was in part supported by a grant no. NSF-RTF-DMS-1840265 from the National Science Foundation, Computational Mathematics Program.

	\pagebreak
	\newpage
	\bibliographystyle{plain}
	\bibliography{main.bib}

\end{document}